\documentclass[aop,MSNbibl,seceqn,citesort,dvips]{arximspdf}
\usepackage{mathrsfs}
\usepackage{eufrak}

\doi{10.1214/10-AOP531}
\volume{38}
\issue{5}
\pubyear{2010}
\firstpage{1947}
\lastpage{1985}

\makeatletter
\newtheorem{prop}[defi]{Proposition}

\newtheorem{theorem}[defi]{Theorem}
\newtheorem{proposition}[defi]{Proposition}
\newtheorem{lemma}[defi]{Lemma}
\newtheorem{cor}[defi]{Corollary}
\newproclaim{defi}{Definition}[section]
\newproclaim{fact}{Fact}
\newproclaim{rem}[defi]{Remark}
\newproclaim{example}[defi]{Example}

\def\HH{\EuFrak H}

\makeatother

\begin{document}
\begin{frontmatter}

\title{Invariance principles for homogeneous sums: Universality of
Gaussian Wiener chaos}
\runtitle{Invariance principles for homogeneous sums}

\begin{aug}
\author[A]{\fnms{Ivan} \snm{Nourdin}\corref{}\ead[label=e1]{ivan.nourdin@upmc.fr}},
\author[B]{\fnms{Giovanni} \snm{Peccati}\ead[label=e2]{giovanni.peccati@gmail.com}}
\and
\author[C]{\fnms{Gesine} \snm{Reinert}\ead[label=e3]{reinert@stats.ox.ac.uk}}
\runauthor{I. Nourdin, G. Peccati and G. Reinert}
\affiliation{Universit\'e Paris VI, Universit\'{e} du Luxembourg and
Oxford University}
\address[A]{I. Nourdin\\
Laboratoire de Probabilit\'es et \\
\quad Mod\`eles Al\'eatoires\\
Universit\'e Pierre et Marie Curie (Paris VI)\\
Bo\^ite courrier 188, 4 place Jussieu\\
75252 Paris Cedex 05\\
France\\
\printead{e1}} %adresu isvedimo komanda gale!
\address[B]{G. Peccati\\
Unit\'{e} de Recherche en Math\'{e}matiques\\
Universit\'{e} du Luxembourg\\
162A, avenue de la Fa\"{\i}encerie\\
L-1511 Luxembourg\\
Grand-Duchy of Luxembourg.\\
On leave from: Universit\'{e} Paris Ouest --\\
\quad Nanterre la D\'{e}fense,
France.\\
\printead{e2}}
\address[C]{G. Reinert\\
Department of Statistics\\
University of Oxford\\
1 South Parks Road\\
Oxford OX1 3TG\\
United Kingdom\\
\printead{e3}}
\end{aug}

% HISTORY:
\received{\smonth{4} \syear{2009}}
\revised{\smonth{1} \syear{2010}}

% ABSTRACT
%
\begin{abstract}
We compute explicit bounds in the normal and chi-square approximations
of multilinear homogenous sums (of arbitrary order) of general centered
independent random variables with unit variance.
In particular, we show that chaotic random variables enjoy the
following form of \textit{universality}: (a)
the normal and chi-square approximations of any homogenous sum can be
completely characterized and assessed by
first switching to its Wiener chaos counterpart, and (b) the simple
upper bounds and convergence criteria
available on the Wiener chaos extend almost verbatim to the class of
homogeneous sums.
\end{abstract}

% KEYWORDS
%

\end{frontmatter}

%s1 ###
\section{Introduction}\label{intr}
%s1.1 ###
\subsection{Overview}
The aim of this paper is to study and characterize the normal and
chi-square approximations of the laws
of \textit{multilinear homogeneous sums} involving general independent
random variables. We shall
perform this task by implicitly combining three probabilistic
techniques, namely: (i) the \textit{Lindeberg
invariance principle} (in a version due to Mossel et al. \cite{MOO}),
%Rotar' \cite{Rotar2} and Mossel \cite{Mossel}),
(ii)
\textit{Stein's method} for the normal and chi-square approximations (see, e.g.,
\cite{ChenShaosur,Reinertsur,Steinorig,Steinbook}), and (iii)~the
\textit{Malliavin calculus of variations} on a Gaussian space (see, e.g.,
\cite{MallBook,Nbook}). Our analysis reveals that the Gaussian Wiener
chaos (see Section \ref{S:WienerC} below for precise definitions)
enjoys the following properties: (a) the normal and chi-square
approximations of any
multilinear homogenous sum are completely characterized and assessed by
those of its Wiener chaos
counterpart, and (b) the strikingly simple upper bounds and
convergence criteria available on the
Wiener chaos (see \cite
{NouPe08,NouPecptrf,NouPeAOP2,NouPeRev,NouPeReinPOIN,NP}) extend almost
verbatim to the class of homogeneous sums. Our findings partially rely
on the notion of ``low influences''
(see again \cite{MOO}) for real-valued functions defined on product
spaces. As indicated by the title, we
regard the two properties (a) and (b) as an instance
of the
\textit{universality phenomenon},
according to which most information about large random systems (such as
the ``distance to Gaussian''
of nonlinear functionals of large samples of independent random
variables) does not depend on the particular distribution
of the components. Other recent examples of the universality phenomenon
appear in the already quoted
paper \cite{MOO}, as well as in the Tao--Vu proof of the circular law
for random matrices, as detailed in
\cite{TV} (see also the Appendix to \cite{TV} by Krishnapur). Observe
that, in Section \ref{S:multivariateext},
we will
prove analogous results for the multivariate normal approximation of
\textit{vectors} of homogenous sums of
possibly different orders. In a further work by the first two authors
(see \cite{NouPeMatrix}) the results of the present paper are applied
in order to deduce universal Gaussian fluctuations for traces
associated with non-Hermitian matrix ensembles.

%s1.2 ###
\subsection{The approach}
In what follows, every random object is defined on a suitable (common)
probability space $(\Omega, \mathscr{F}, P)$. The symbol $E$ denotes
expectation with respect to $P$. We start by giving a precise
definition of the main objects of our study.
\begin{defi}[(Homogeneous sums)]\label{Def:HomSums}
Fix some integers $N,d\geq2$ and write $[N]=\{1,\ldots,N\}$. Let
$\mathbf{X}
= \{X_i\dvtx i\geq1\} $ be a collection of centered independent random
variables, and let $f\dvtx [N]^d \rightarrow\mathbb{R}$ be a \textit{symmetric
function vanishing on diagonals} [i.e., $f(i_1,\ldots,i_d)=0$ whenever
there exist $k\neq j$ such that $i_k=i_j$]. The random variable
%
%
%e2 ###
%e1 ###
\begin{eqnarray}\label{EQ:HOMsums}
 Q_d(N,f,\mathbf{X})&=& Q_d(\mathbf{X})= \sum_{1\leq i_1,\ldots,i_d
\leq
N}f(i_1,\ldots,i_d)X_{i_1}\cdots X_{i_d} \nonumber\\
&=& d! \sum_{\{i_1,\ldots,i_d\}\subset[N]^d}
f(i_1,\ldots,i_d)X_{i_1}\cdots X_{i_d}\\
&=& d! \sum_{1\leq i_1 <\cdots<i_d\leq N}
f(i_1,\ldots,i_d)X_{i_1}\cdots X_{i_d}\nonumber
\end{eqnarray}
is called the \textit{multilinear homogeneous sum}, of order $d$, based on
$f$ and on the first~$N$
elements of $\mathbf{X}$.
\end{defi}

As in (\ref{EQ:HOMsums}), and when there is no risk of confusion, we
will drop the dependence on $N$ and $f$ in order to simplify the notation.
Plainly, $E[Q_d(\mathbf{X})]=0$ and also, if $E(X_i^2)=1$ for every
$i$, then
$
E[Q_d(\mathbf{X})^2] = d!\|f\|^2_{d},
$
where we use the notation
$
\|f\|^2_{d} = \sum_{1\leq i_1,\ldots,i_d \leq N}f^2(i_1,\ldots,i_d)
$
(here and for the rest of the paper). In the following, we will
systematically use the expression ``homogeneous sum'' instead of
``multilinear homogeneous sum.''

Objects such as (\ref{EQ:HOMsums}) are sometimes called ``polynomial
chaoses,'' and play a
central role in several branches of probability theory and stochastic
analysis. When $d=2$, they are typical
examples of quadratic forms. For general $d$, homogeneous sums are, for
example, the basic building blocks of the Wiener, Poisson and Walsh chaoses
(see, e.g., \cite{Privault}). Despite the almost ubiquitous nature of
homogeneous sums,
results concerning the normal approximation of quantities such as (\ref
{EQ:HOMsums}) in the nonquadratic case
(i.e., when $d\geq3$) are surprisingly scarce: indeed, to our
knowledge, the only general statements
in this respect are contained in references~\cite
{deJongBook,deJongMulti}, both by P. de Jong
(as discussed below), and
in a different direction, general criteria allowing to assess the
proximity of the laws of homogenous sums
based on different independent
sequences are obtained in \cite{MOO,Rotar1,Rotar2}.

In this paper we are interested in controlling objects of the
type\break
$d_{\mathscr H}\{Q_d(\mathbf{X}); Z\}$,
where: (i) $Q_d(\mathbf{X})$ is defined in (\ref{Def:HomSums}), (ii) $Z$
is either a standard Gaussian
$\mathscr{N}(0,1)$ or a centered chi-square random variable, and (iii)
the distance $d_{\mathscr{H}}\{F;G\}$,
between the laws of two random variables $F$ and $G$, is given by
%
%
%e3 ###
\begin{equation}\label{EQ:distance}
d_{\mathscr{H}}\{F; G\}=\sup\{| E[h(F)] - E[h(G)]
|\dvtx h\in\mathscr{H}\}
\end{equation}
with $\mathscr{H}$ some suitable class of real-valued functions. Even
with some uniform control on the components
of $\mathbf{X}$, the problem of directly and generally assessing
$d_{\mathscr H}\{Q_d(\mathbf{X}); Z\}$ looks very
arduous. Indeed, any estimate comparing the laws of $Q_d(\mathbf{X})$ and
$Z$ capriciously depends on the kernel
$f$, and on the way in which the analytic structure of $f$ interacts
with the specific ``shape'' of the
distribution of the random variables $X_i$. One revealing picture of
this situation appears if one tries to
evaluate the moments of $Q_d(\mathbf{X})$ and to compare them with those
of $Z$; see, for example, \cite{PecTaqSURV}
for a discussion of some associated combinatorial structures. In the
specific case where $Z$ is Gaussian, one should also observe that
$Q_d(\mathbf{X})$ is a \textit{completely degenerate $U$-statistic}, as
$E[f(i_1,\ldots,i_d) X_{i_1} x_{i_2} \cdots x_{i_d} ] =0$ for all
$x_{i_2}, \ldots, x_{i_d}$, so that the standard results for the normal
approximation of $U$-statistics do not apply.

The main point developed in the present paper is that one can
successfully overcome these difficulties by implementing
the following strategy: first (I) measure the distance
$d_{\mathscr H}\{Q_d(\mathbf{X});
Q_d(\mathbf{G})\},$ between the law of $Q_d(\mathbf{X})$ and the law of
the random variable $Q_d(\mathbf{G})$, obtained
by replacing $\mathbf{X}$ with a centered standard i.i.d. Gaussian sequence
$\mathbf{G} = \{G_i \dvtx i\geq1\}$; then (II)
assess the distance $d_{\mathscr H}\{Q_d(\mathbf{G}); Z\}$; and
finally (III) use the triangle
inequality in order to write
%
%
%e4 ###
\begin{equation}\label{trianglebound}
d_{\mathscr H}\{Q_d(\mathbf{X}); Z\}\leq d_{\mathscr H}\{
Q_d(\mathbf{X}); Q_d(\mathbf{G})\}+ d_{\mathscr H}\{Q_d(\mathbf{G}); Z\}.
\end{equation}
We will see in the subsequent sections that the power of this approach
resides in the following two facts.

\begin{fact} The distance evoked at Point (I)
can be
effectively controlled by means of the techniques
developed in \cite{MOO}, where the authors have produced a general theory
allowing to estimate the
distance between homogeneous sums constructed from different sequences
of independent random variables. A full discussion of this point is
presented in Section \ref{S:MOO} below. In Theorem \ref{thmMOO} we
shall observe that, under the assumptions that $E(X_i^2)=1$ and that
the moments $E(|X_i|^{3})$ are uniformly bounded by some constant
$\beta>0$ (recall that the $X_i$'s are centered), one can deduce from
\cite{MOO} (provided that the elements of $\mathscr{H}$ are
sufficiently smooth) that
%
%
%e5 ###
\begin{equation}\label{Eq:moo1}
\qquad d_{\mathscr H}\{Q_d(\mathbf{X}); Q_d(\mathbf{G})\} \leq C \times
\sqrt
{\max_{1\leq i\leq N}
\sum_{\{i_2,\ldots,i_{d}\}\in[N]^{d-1} }f^2(i,i_2,\ldots,i_{d})},
\end{equation}
where $C$ is a constant depending only on $d$, $\beta$ and on the class
$\mathscr{H}$.
The quantity
%
%
%e6 ###
\begin{eqnarray}\label{infi}
\operatorname{Inf}_i(f)&:=&\sum_{\{i_2,\ldots,i_{d}\}\in[N]^{d-1}
}f^2(i,i_2,\ldots,i_{d})
\nonumber
\\[-8pt]
\\[-8pt]
\nonumber
&\phantom{:}=&\frac{1}{(d-1)!} \sum_{1\leq
i_2,\ldots,i_{d}\leq N
}f^2(i,i_2,\ldots,i_{d})
\end{eqnarray}
is called the \textit{influence} of the variable $i$, and roughly
quantifies the contribution of~$X_i$ to the overall configuration of
the homogenous sum $Q_d(\mathbf{X})$. Influence indices already appear
(under a different name) in the papers by Rotar' \cite{Rotar1,Rotar2}.
\end{fact}

\begin{fact} The random variable $Q_d(\mathbf{G})$ is an
element of the $d$th \textit{Wiener chaos} associated with $\mathbf{G}$
(see Section \ref{S:WienerC} for definitions). As such, the distance
between $Q_d(\mathbf{G})$ and $Z$ (in both the normal and the chi-square
cases) can be assessed by means of the results appearing in
\cite{NouPe08,NouPecptrf,NouPeAOP2,NouPeRev,NouPeReinPOIN,NO,NP,PT},
which are in turn based on a
powerful interaction between standard Gaussian analysis, Stein's method
and the Malliavin calculus on
variations. As an example, Theorem \ref{T:4thcumulant} of Section
\ref{S:newbounds} proves that,
if $Q_d(\mathbf{G})$ has variance one and $Z$ is standard Gaussian, then
%
%
%e7 ###
\begin{equation}\label{bound1}
\qquad d_{\mathscr H}\{Q_d(\mathbf{G}); Z\} \leq C\sqrt{|E[Q_d(\mathbf{
G})^4]-E(Z^4)|}=C\sqrt{|E[Q_d(\mathbf{G})^4]-3|} ,
\end{equation}
where $C>0$ is some finite constant depending only on $\mathscr{H}$
and $d$.
\end{fact}
%Moreover, if the elements of $\mathbf{X}$ have bounded fourth moments,
%then the RHS of (\ref{bound1}) can be directly estimated by means of
%the fourth moment of $Q_d(\mathbf{X})$ and some maxima over the
%influence
%functions appearing in (\ref{Eq:moo1}) (see Section

%s1.3 ###
\subsection{Universality}\label{univers}

Bounds such as (\ref{Eq:moo1}) and (\ref{bound1}) only partially
account for the term ``universality''
appearing in the title of the present paper. Our techniques allow
indeed to prove the following statement,
involving vectors of homogeneous sums of possibly different orders; see
also Theorem
\ref{T:UNIVmultivbis} for a more general statement.

\begin{theorem}[(Universality of Wiener chaos)]\label{T:UNIVmultiv}
Let $\mathbf{G}=\{G_i \dvtx  i\geq1\}$ be a
standard centered i.i.d. Gaussian sequence, and fix integers
$m\geq1$
and $d_1,\ldots,d_m\geq2$.
For every $j=1,\ldots,m$, let $\{(N^{(j)}_n, f^{(j)}_n) \dvtx  n\geq1\}$
be a
sequence such that $\{N^{(j)}_n \dvtx  n\geq1\}$
is a sequence of integers going to infinity, and each function
$f_n^{(j)}\dvtx  [N_n^{(j)}]^{d_j} \rightarrow\mathbb{R}$
is symmetric and vanishes on diagonals.
Define $Q_{d_j}(N_n^{(j)}, f^{(j)}_n, \mathbf{G})$, $n\geq1$,
according to
\textup{(\ref{EQ:HOMsums})}
and assume that,
for every $j=1,\ldots,m$, the sequence $E[Q_{d_j}(N_n^{(j)},f^{(j)}_n, \mathbf{
G})^2]$, $n\geq1$, is bounded.
Let $V$ be a $m\times m$ nonnegative symmetric matrix whose diagonal
elements are different from zero, and let $\mathscr{N}_m(0,V)$
indicate a centered Gaussian vector with covariance $V$. Then, as $n
\rightarrow\infty$, the following conditions (1) and (2)
are equivalent: (1)~The vector $\{Q_{d_j}(N_n^{(j)},f^{(j)}_n,
\mathbf{G}) \dvtx  j=1,\ldots,m\}$ converges in law to $\mathscr{N}_m(0,V)$;
(2)~for every sequence $\mathbf{X } = \{X_i \dvtx  i\geq1\}$ of independent
centered random variables,
with unit variance and such that $\sup_i E|X_i|^{3} <\infty$, the law
of the vector
$\{Q_{d_j}(N_n^{(j)},f^{(j)}_n, \mathbf{X}) \dvtx  j=1,\ldots,m\}$ converges to
the law of $\mathscr{N}_m(0,V)$ in the Kolmogorov distance.
\end{theorem}

\begin{rem} 1. Given random vectors $F = (F_1,\ldots,F_m)$ and $H =
(H_1,\ldots,\break H_m)$, $m\geq1$, the \textit{Kolmogorov distance} between the
law of $F$ and the law of $H$ is defined as
%
%
%e8 ###
\begin{eqnarray}\label{IntroKOLdis}
d_{\mathrm{Kol}}(F,H) &=& \sup_{(z_1,\ldots,z_m)\in\mathbb{R}^m} |
P(F_1\leq
z_1,\ldots,F_m \leq z_m)
\nonumber
\\[-8pt]
\\[-8pt]
\nonumber
&&\hspace*{35pt}\qquad{}- P(H_1\leq z_1,\ldots,H_m \leq
z_m) |.
\end{eqnarray}
Recall that the topology induced by $d_{\mathrm{Kol}}$ on the class of all
probability measures on $\mathbb{R}^m$ is strictly stronger than the topology
of convergence in distribution.\vspace*{-6pt}
\begin{longlist}
\item[2.] Note that, in the statement of Theorem \ref{T:UNIVmultiv},
we do not require that the matrix
$V$ is \textit{positively definite}, and we do \textit{not} introduce
any assumption on the asymptotic behavior of influence indices.
\item[3.] Due to the matching moments up to second order, one has that
\begin{eqnarray*}
&&E\bigl[Q_{d_i}\bigl(N_n^{(i)},f^{(i)}_n, \mathbf{G}\bigr) \times
Q_{d_j}\bigl(N_n^{(j)},f^{(j)}_n, \mathbf{G}\bigr)\bigr]\\
&&\qquad= E\bigl[Q_{d_i}\bigl(N_n^{(i)},f^{(i)}_n, \mathbf{X}\bigr) \times
Q_{d_j}\bigl(N_n^{(j)},f^{(j)}_n, \mathbf{X}\bigr)\bigr]
\end{eqnarray*}
for every $i,j=1,\ldots,m$ and every sequence $\mathbf{X}$ as in Theorem
\ref{T:UNIVmultiv}.
%given in Theorem \ref{T:byebyedeJong} below.
\end{longlist}
\end{rem}

Theorem \ref{T:UNIVmultiv} basically ensures that any statement
concerning the asymptotic normality of
(vectors of) general homogeneous sums can be proved by simply focusing
on the elements of a Gaussian Wiener
chaos. Since central limit theorems (CLTs) on Wiener chaos are by now
completely characterized (thanks to the
results proved in \cite{NP}), this fact represents a clear methodological
breakthrough. As explained later in the paper, and up to the
restriction on the third moments,
we regard Theorem \ref{T:UNIVmultiv} as the first exact
equivalent---for homogeneous sums---of
the usual CLT for linear functionals of i.i.d. sequences. The proof of
Theorem \ref{T:UNIVmultiv} is
achieved in Section \ref{S:multivariateext}.

\begin{rem}
When dealing with the multidimensional case, our way to use the
techniques developed in
\cite{Mossel} makes it
unavoidable to require a uniform bound on the third moments of $\mathbf{
X}$. However, one advantage
is that we easily obtain convergence in the Kolmogorov distance, as
well as explicit
upper bounds on the rates of convergence. We will see below (see
Theorem \ref{T:byebyedeJong} for a
precise statement) that in the one-dimensional case one can simply
require a bound on the moments of order
$2+\varepsilon$, for some $\varepsilon>0$. Moreover, still in the
one-dimensional case and when the sequence
$\mathbf{X}$ is i.i.d., one can alternatively deduce convergence in
distribution from a result by
Rotar' (\cite{Rotar2}, Proposition 1), for which the existence of moments
of order greater than $2$ is not
required.
\end{rem}

%s1.4 ###
\subsection{The role of contractions}\label{SS:role}
The universality principle stated in Theorem \ref{T:UNIVmultiv} is
based on \cite{MOO}, as well as
on general characterizations of (possibly multidimensional) CLTs on a
fixed Wiener chaos. Results of
this kind have been first proved in \cite{NP} (for the one-dimensional
case) and \cite{PT}
(for the multidimensional case), and make an important use of the
notion of ``contraction''
of a given deterministic kernel. When studying homogeneous sums, one is
naturally led to deal with contractions defined on discrete sets of the
type $[N]^d$, $N\geq1$.
In this section we shall briefly explore these discrete objects, in
particular, by pointing out that discrete contractions are indeed the
key element in the proof of Theorem \ref{T:UNIVmultiv}. More general
statements, as well as complete proofs, are given in Section \ref{S:newbounds}.

\begin{defi}\label{def-contr}
Fix $d,N\geq2$. Let $f \dvtx  [N]^d \rightarrow\mathbb{R}$ be a
symmetric
function vanishing of diagonals.
For every $r= 0,\ldots,d$, the \textit{contraction} $f\star_r f$ is the
function on $[N]^{2d-2r}$ given by
\begin{eqnarray*}
&&f \star_r f(j_1,\ldots,j_{2d-2r}) \\
&&\qquad=
{\sum_{1\leq a_1,\ldots,a_r \leq N}}
f(a_1,\ldots,a_r,j_1,\ldots,j_{d-r})
f(a_1,\ldots,a_r,j_{d-r+1},\ldots,j_{2d-2r}).
\end{eqnarray*}
\end{defi}
Observe that $f\star_r f $ is not necessarily symmetric and does not
necessarily vanish on diagonals. The symmetrization of $f\star_r f$ is
written $f\,\widetilde{\star}_r\, f$.
The following result, whose proof is achieved in Section \ref
{S:multivariateext} as a special case
of Theorem \ref{T:UNIVmultivbis}, is based on the findings of \cite{NP,PT}.

\begin{prop}[(CLT for chaotic sums)]\label{T:NGCGHomSums} Let the
assumptions and notation of Theorem \textup{\ref{T:UNIVmultiv}} prevail, and
suppose, moreover, that, for every $i,j=1,\ldots,m$ (as $n\rightarrow
\infty$),
%
%
%e9 ###
\begin{equation}\label{covcv}
E\bigl[Q_{d_i}\bigl(N_n^{(i)},f^{(i)}_n, \mathbf{G}\bigr)\times
Q_{d_j}\bigl(N_n^{(j)},f^{(j)}_n, \mathbf{G}\bigr)\bigr] \rightarrow V(i,j),
\end{equation}
where $V$ is a nonnegative symmetric matrix. Then, the following three
conditions \textup{(1)--(3)} are equivalent, as $n\rightarrow\infty$:
\textup{(1)} The vector $\{Q_{d_j}(N_n^{(j)},f^{(j)}_n, \mathbf{G}) \dvtx
j=1,\ldots,m\}$ converges in law to a centered Gaussian vector with
covariance matrix $V$;
\textup{(2)} for every $j=1,\ldots,m$, $E[Q_{d_j}(N_n^{(j)},f^{(j)}_n,
\mathbf{G})^4] \rightarrow3V(i,i)^2$;
\textup{(3)} for every $j=1,\ldots,m$ and every $r=1,\ldots,d_j -1$, $\|
f_n^{(j)}\star_r f_n^{(j)}\|_{2d_j-2r} \rightarrow0$.
\end{prop}

\begin{rem}
Strictly speaking, the results of \cite{PT} only deal with the case
where $V$ is positive definite.
The needed general result will be obtained in Section~\ref
{S:multivariateext} by means of Malliavin
calculus.
\end{rem}

Let us now briefly sketch the proof of Theorem \ref{T:UNIVmultiv}.
Suppose that the sequence $E[Q_{d_j}(N_n^{(j)},f^{(j)}_n, \mathbf{G})^2]$
is bounded and that
the vector $\{Q_{d_j}(N_n^{(j)},f^{(j)}_n,\break \mathbf{G}) \dvtx  j=1,\ldots,m\}$
converges in law to $\mathscr{N}_m(0,V)$.
Then, by uniform integrability (using Proposition \ref{P:Hyper}), the
convergence (\ref{covcv}) is satisfied and, according to Proposition~\ref{T:NGCGHomSums},
we have $\|f_n^{(j)}\star_{d_j-1} f_n^{(j)}\|_{2} \rightarrow0$. The
crucial remark is now that
%
%
%e11 ###
%e10 ###
\begin{eqnarray}\label{EQ:cruxcontr}
\bigl\|f_n^{(j)}\star_{d_j-1} f_n^{(j)}\bigr\|^2_{2} &\geq&\sum
_{1\leq i\leq
N_n^{(j)}} \biggl[
\sum_{1\leq i_2,\ldots,i_{d_j}\leq N_n^{(j)}}
f_n^{(j)}(i,i_2,\ldots,i_{d_j})^2\biggr]^2 \nonumber\\
%&\geq& \sum_{i=1}^{N_n^{(j)}} [
%f_n^{(j)}(i,i_2,\ldots,i_{d_j})^2]^2 \\
&\geq& \max_{1\leq i\leq N_n^{(j)}} \biggl[
\sum_{1\leq i_2,\ldots,i_{d_j}\leq N_n^{(j)}}
f_n^{(j)}(i,i_2,\ldots,i_{d_j})^2\biggr]^2\\
&=& \Bigl[(d_j-1)! \max_{1\leq i\leq N_n^{(j)}}
\operatorname{Inf}_i\bigl(f_n^{(j)}\bigr)\Bigr]^2\nonumber
\end{eqnarray}
[recall formula (\ref{infi})], from which one immediately obtains that,
as $n\rightarrow\infty$,
%
%
%e12 ###
\begin{equation}\label{multiINF}
\max_{1\leq i\leq N_n^{(j)}} \operatorname{Inf}_i\bigl(f_n^{(j)}\bigr)
\rightarrow0\qquad
\mbox{for every } j=1,\ldots,m.
\end{equation}
The proof of Theorem \ref{T:UNIVmultiv} is concluded by using Theorem
\ref{T:MULTIMOO},
which is a statement in the same vein as the results established in
\cite{Mossel}, that is, a
multidimensional version of the findings of \cite{MOO}. Indeed, this
result will imply that, if
(\ref{multiINF}) is verified, then, for every sequence $\mathbf{X}$ as in
Theorem \ref{T:UNIVmultiv}, the distance between the law of $\{
Q_{d_j}(N_n^{(j)},f^{(j)}_n, \mathbf{G}) \dvtx
j=1,\ldots,m\}$ and the law of $\{Q_{d_j}(N_n^{(j)},f^{(j)}_n, \mathbf{X}) \dvtx
j=1,\ldots,m\}$ necessarily tends to
zero and, therefore, the two sequences must converge in distribution to
the same limit.

As proved in \cite{NouPe08}, contractions play an equally important
role in the chi-square
approximation of the laws of elements of a fixed chaos of even order.
Recall that a random variable
$Z_\nu$ has a \textit{centered chi-square distribution}\vspace*{-1pt} with $\nu
\geq1$
degrees of freedom
[noted $Z_\nu\sim\chi^2(\nu)$] if $Z_\nu\stackrel{\mathrm{Law}}{=}
\sum_{i=1}^\nu(G_i^2 -1)$,
where $(G_1,\ldots,G_\nu)$ is a vector of i.i.d. $\mathscr{N}(0,1)$ random
variables. Note that
$E(Z^2_\nu)=2\nu$, $E(Z^3_\nu)=8\nu$ and $E(Z^4_\nu)=12\nu
^2+48\nu$.

\begin{theorem}[(Chi-square limit theorem for chaotic sums, \cite
{NouPe08})]\label{T:IvGioNCLT}
Let $\mathbf{G}=\{G_i \dvtx  i\geq1\}$ be
a standard centered i.i.d. Gaussian sequence, and fix an even integer
$d\geq2$.
Let $\{N_n, f_n \dvtx  n\geq1\}$ be a sequence such that $\{N_n \dvtx
n\geq1\}
$ is a sequence of integers going to
infinity, and each $f_n\dvtx  [N_n]^{d} \rightarrow\mathbb{R}$ is
symmetric and
vanishes on diagonals.
Define $Q_{d}(N_n,f_n, \mathbf{G})$, $n\geq1$, according to \textup{(\ref
{EQ:HOMsums})}, and assume that, as
$n\rightarrow\infty$, $E[Q_{d}(N_n,f_n, \mathbf{G})^2]\rightarrow2\nu$.\vspace*{-1pt}
Then, as $n\rightarrow\infty$, the following conditions \textup{(1)--(3)} are equivalent:
\textup{(1)} $Q_{d}(N_n,f_n, \mathbf{G})\stackrel{\mathrm{Law}}{\rightarrow}
Z_\nu\sim\chi^2(\nu)$;
\textup{(2)}~$E[Q_{d}(N_n,f_n, \mathbf{G})^4]-12E[Q_{d}(N_n,f_n, \mathbf{
G})^3]\rightarrow E[Z_\nu^4]-12E[Z_\nu^3] = 12\nu^2-48\nu$;
\textup{(3)} $\|f_n\,\widetilde{\star}_{d/2}\, f_n-c_d\times f_n\|_d\to0$
and $\|f_n\star_{r} f_n\|_{2d-2r}\rightarrow0$ for every
$r=1,\ldots,d-1$ such that $r\neq d/2$, where
$
c_{d}:=
4(d/2)!^3d!^{-2}.
$
\end{theorem}

%s1.5 ###
\subsection{Example: Revisiting de Jong's criterion}\label{sec-dejong}
To further clarify the previous discussion, we provide an illustration
of how one can use our results in order to refine a remarkable result
by de Jong, originally proved in \cite{deJongMulti}.

\begin{theorem}[(See \cite{deJongMulti})]\label{T:deJong}
Let $\mathbf{X } = \{X_i \dvtx  i\geq1 \}$ be a sequence of independent
centered random variables
such that $E(X_i^2)=1$ and $E(X_i^4) < \infty$ for every~$i$.
Fix $d\geq2$, and let $\{N_n, f_n \dvtx  n\geq1\}$ be a sequence
such that
$\{N_n \dvtx  n\geq1\}$ is a sequence of integers going to
infinity, and each $f_n\dvtx  [N_n]^{d} \rightarrow\mathbb{R}$ is
symmetric and
vanishes on diagonals.
Define $Q_d(n,\mathbf{X })=Q_{d}(N_n,f_n, \mathbf{X})$, $n\geq1$, according
to~\textup{(\ref{EQ:HOMsums})}.
Assume that $E[Q_d(n,\mathbf{X })^2]=1$ for all $n$. Suppose that, as
$n\rightarrow\infty$: \textup{(i)}~$E[Q_d(n,\mathbf{X })^4]\rightarrow3$, and
\textup{(ii)} $\max_{1\leq i\leq N_n} \operatorname{Inf}_i(f_n) \rightarrow
0$. Then,
$Q_d(n,\mathbf{X })$ converges in law to
$Z\sim\mathscr{N}(0,1)$.
\end{theorem}

%Theorem \ref{T:deJong} extends to homogeneous sums that
%are not necessarily multilinear. Also, the fact that $E(X^4_i)<\infty$
%does not appear in \cite{deJongMulti} (however, an inspection of the
%proof of Theorem \ref{T:deJong} shows that this assumption is indeed
%necessary).}

In the original proof given in \cite{deJongMulti}, assumption (i) in
Theorem \ref{T:deJong} appears as a
convenient (and mysterious) way of reexpressing the asymptotic ``lack
of interaction'' between products of the
type $X_{i_1}\cdots X_{i_d}$, whereas assumption (ii) plays
the role of a usual Lindeberg-type
assumption. In the present paper, under the slightly stronger
assumption that $\sup_i E(X_i^4) < \infty$,
we will be able to produce bounds neatly indicating the exact roles of
both assumptions (i) and (ii). To see this, define~$d_{\mathscr{H}}$
according to (\ref{EQ:distance}), and set $\mathscr{H}$ to be the
class of
thrice differentiable functions whose first three derivatives are
bounded by some finite constant $B>0$.
In Section~\ref{S:normalHOMSUMS}, in the proof of Theorem \ref
{mickey}, we will show that there exist
universal, explicit, finite constants $C_1,C_2,C_3>0$, depending only
on $\beta$, $d$ and $B$, such that (writing $\mathbf{G}$ for an i.i.d.
centered standard Gaussian sequence)\looseness=1
%
%
%e15 ###
%e14 ###
%e13 ###
\begin{eqnarray}
d_{\mathscr{H}}\{Q_d(n,\mathbf{X}); Q_d(n,\mathbf{G})\}& \leq& C_1
\times
\sqrt{\max_{1\leq i\leq N_n} \operatorname{Inf}_i(f_n)} \label
{J1},\\
d_{\mathscr{H}}\{Q_d(n,\mathbf{G}); Z \} & \leq& C_2 \times\sqrt
{|E[Q_d(n,\mathbf{G})^4] - 3|}, \label{J2}\\
\quad |E[Q_d(n,\mathbf{X})^4] - E[Q_d(n,\mathbf{G})^4]| & \leq& C_3
\times
\sqrt{\max_{1\leq i\leq N_n} \operatorname{Inf}_i(f_n)}. \label{J3}
\end{eqnarray}
In particular, the estimates (\ref{J1}) and (\ref{J3}) show that
assumption (ii) in Theorem~\ref{T:deJong}
ensures that both the laws and the fourth moments of $Q_d(n,\mathbf{X})$
and $Q_d(n,\mathbf{G})$ are asymptotically
close: this fact, combined with assumption (i), implies that the LHS of
(\ref{J2}) converges to zero, hence
so does $d_{\mathscr{H}}\{Q_d(n,\mathbf{X}); Z \}$. This gives
an alternate proof of Theorem \ref{T:deJong} in the case of uniformly
bounded fourth moments.

Also, by combining the universality principle stated in Theorem \ref
{T:UNIVmultiv} with (\ref{J2})
(or, alternatively, with Proposition \ref{T:NGCGHomSums} in the case
$m=1$), one obtains the following ``universal version'' of de Jong's criterion.

\begin{theorem}\label{T:byebyedeJong}
Let $\mathbf{G } = \{X_i \dvtx  i\geq1 \}$ be a centered i.i.d. Gaussian
sequence with unit variance.
Fix $d\geq2$, and let $\{N_n, f_n \dvtx  n\geq1\}$ be a sequence
such that
$\{N_n \dvtx  n\geq1\}$ is a sequence of integers going to
infinity, and each $f_n\dvtx  [N_n]^{d} \rightarrow\mathbb{R}$ is
symmetric and
vanishes on diagonals.
Define $Q_d(n,\mathbf{G})=Q_{d}(N_n,f_n, \mathbf{G})$, $n\geq1$, according
to \textup{(\ref{EQ:HOMsums})}.
Assume that $E[Q_d(n,\mathbf{G})^2]\to1$ as $n\to\infty$.
Then, the following four properties are equivalent as $n\rightarrow
\infty$:
\begin{enumerate}[(1)]
\item[(1)] The sequence $Q_d(n,\mathbf{G})$ converges in law to
$Z\sim
\mathscr{N}(0,1)$.
\item[(2)] $E[Q_d(n,\mathbf{G})^4]\rightarrow3$.
\item[(3)] For every sequence $\mathbf{X } = \{X_i \dvtx  i\geq1\}$ of
independent centered random variables
with unit variance and such that $\sup_i E|X_i|^{2+\varepsilon} <\infty$
for some $\varepsilon>0$, the sequence $Q_d(n,\mathbf{X})$ converges in law
to $Z\sim\mathscr{N}(0,1)$ in the Kolmogorov distance.
\item[(4)] For every sequence $\mathbf{X } = \{X_i \dvtx  i\geq1\}$ of
independent and identically distributed centered random variables
with unit variance, the sequence\break $Q_d(n,\mathbf{X})$ converges in law to
$Z\sim\mathscr{N}(0,1)$ (not necessarily in the Kolmogorov distance).
\end{enumerate}
\end{theorem}

\begin{rem}
1. Note that at point (4) of the above statement we \textit{do not}
require the existence of moments of
order greater than 2. We will see that the equivalence between (1) and
(4) is partly a consequence of Rotar's
results (see~\cite{Rotar2}, Proposition 1).\vspace*{-6pt}
\begin{longlist}
\item[2.] Theorem \ref{T:byebyedeJong} is a particular case of
Theorem \ref{T:UNIVmultiv}, and can be seen as refinement of de
Jong's Theorem \ref{T:deJong}, in the sense that: (i) since several
combinatorial devices are at hand (see, e.g., \cite{PecTaqSURV}), it
is in general easier to evaluate moments of multilinear forms of
Gaussian sequences than of general sequences, and (ii)~when the $\{X_i\}
$ are not identically distributed, we only need existence (and uniform
boundedness) of the moments of order $2+\varepsilon$.
\end{longlist}
\end{rem}

In Section \ref{S:multivariateext} we will generalize the content of
this section to multivariate Gaussian approximations. By using
Proposition \ref{T:IvGioNCLT} and \cite{Rotar2}, Proposition 1, one
can also obtain the following universal chi-square limit result.

\begin{theorem}\label{T:seeyadeJong} We let the notation of Theorem
\textup{\ref{T:byebyedeJong}} prevail, except that
we now assume that $d\geq2$ is an even integer and $E[Q_d(n,\mathbf{
G})^2] \rightarrow2\nu$, where $\nu\geq1$ is
an integer. Then, the following four conditions \textup{(1)--(4)}
are equivalent as $n\rightarrow\infty$:
\textup{(1)} The sequence $Q_d(n,\mathbf{G})$ converges in law to $Z_\nu\sim
\chi^2(\nu)$;
\textup{(2)}~$E[Q_d(n,\mathbf{G})^4]-12E[Q_d(n,\mathbf{G})^3]\rightarrow
E(Z_\nu
^4)-12E(Z_\nu^3)=12\nu^2-48\nu$;
\textup{(3)}~for every sequence $\mathbf{X } = \{X_i \dvtx  i\geq1\}$ of
independent centered random variables
with unit variance and
such that $\sup_i E|X_i|^{2+\varepsilon} <\infty$ for some $\varepsilon>0$,
the sequence $Q_d(n,\mathbf{X})$ converges in law to $Z_\nu$;
\textup{(4)} for every sequence $\mathbf{X } = \{X_i \dvtx  i\geq1\}$ of
independent and identically distributed centered random variables
with unit variance, the sequence $Q_d(n,\mathbf{X})$ converges in law to
$Z_\nu$.
\end{theorem}

%s1.6 ###
\subsection{Two counterexamples}

\subsubsection*{``There is no universality for sums of order one''} One
striking feature of Theorems \ref{T:UNIVmultiv} and \ref
{T:byebyedeJong} is that
\textit{they do not have any equivalent for sums of order} $d=1$.
To see this, consider an array of real numbers $\{f_n(i) \dvtx  1\leq
i\leq
n\}$ such that
$\sum_{i=1}^n f_n^2(i) =1$. Let $\mathbf{G} = \{G_i \dvtx  i\geq1\}$ and
$\mathbf{
X} = \{X_i \dvtx  i\geq1\}$ be,
respectively, a centered i.i.d. Gaussian sequence with unit variance,
and a sequence of independent random
variables with zero mean and unit variance. Then, $Q_1(n,\mathbf
{G}):=\sum
_{i=1}^n f_n(i)G_i \sim
\mathscr{N}(0,1)$ for every $n$, but it is in general \textit{not true}
that $Q_1(n,\mathbf{X}):=\sum_{i=1}^n
f_n(i)X_i$ converges in law to a Gaussian random variable [just take
$X_1$ to be non-Gaussian, $f_n(1) =1$
and $f_n(j)=0$ for $j>1$]. As it is well known, to ensure that
$Q_1(n,\mathbf{X})$ has a Gaussian limit, one customarily adds the
Lindeberg-type requirement
that $\max_{1\leq i\leq n}| f_n(i)| \rightarrow0$.
%(see e.g. \cite[Th. 4.12]{KallBook}).
A closer inspection indicates
that the fact that no Lindeberg conditions are required in Theorems \ref
{T:UNIVmultiv} and \ref{T:byebyedeJong} is due to the implication (1)
$\Rightarrow$ (3) in Proposition \ref{T:NGCGHomSums},
as well as to the inequality (\ref{EQ:cruxcontr}).

\subsubsection*{``Walsh chaos is not universal''}
One cannot replace the
Gaussian sequence~$\mathbf{G}$ with a Rademacher one in the
statements of Theorems \ref{T:UNIVmultiv} and \ref
{T:byebyedeJong}. Let
$\mathbf{X} = \{X_i \dvtx  i\geq1\}$ be an i.i.d. Rademacher sequence,
and fix
$d\geq2$. For every $N\geq d$, consider the homogeneous sum
$
Q_d( N, \mathbf{X }) =\break X_1X_2\cdots X_{d-1}\sum_{i=d}^N \frac
{X_i}{\sqrt{N-d+1} }.
$
It is easily seen that each $Q_d(N,\mathbf{X})$ can be written in the form
(\ref{EQ:HOMsums}),
for some symmetric $f=f_N$ vanishing on diagonals and such that $d!\|
f_N\|^2_d=1$.
Since $X_1X_2\cdots X_{d-1}$ is a random sign independent of
$\{X_i \dvtx  i\geq d\}$, a
simple application of the central limit theorem yields that, as
$N\rightarrow\infty$, $Q_d(N,\mathbf{X})
\stackrel{\mathrm{Law}}{\rightarrow} \mathscr{N}(0,1)$. On the other hand,
if $\mathbf{G} = \{G_i \dvtx  i\geq1\}$ is a i.i.d. standard Gaussian sequence,
one sees that $Q_d( N, \mathbf{G }) \stackrel{\mathrm{Law}}{=} G_1\cdots G_d$, for every $N\geq2$.
Since (for $d\geq2$) the random variable $G_1\cdots
G_d$ is
not Gaussian, this yields that
$Q_d( N, \mathbf{G })\stackrel{\mathrm{Law}}{\not\rightarrow} \mathscr{N}(0,1)$
as $n\to\infty$.

\begin{rem}
1. In order to enhance the readability of the forthcoming material,
we decided not to state some of our findings in full generality. In
particular: (i)
It will be clear later on that the results of this paper easily extend
to the case of \textit{infinite} homogeneous sums [obtained by putting
$N=+\infty$ in (\ref{EQ:HOMsums})]. This requires, however, a
somewhat heavier notation, as well as some distracting digressions
about convergence. (ii) Our findings do not hinge at all on the fact
that $\mathbb{N}$ is an ordered set: it follows that our results
exactly apply to homogeneous sums of random variables indexed by a
general finite set.\vspace*{-6pt}
%(iii) Our results on chi-square approximations could be slightly
%modified in order to accommodate
%approximations by arbitrary centered Gamma laws. To do this, one can
%use the results about Gamma approximations proved in
%
\begin{longlist}[2.]
\item[2.] As discussed below, the results of this paper are tightly related
with a series of recent findings concerning the normal and Gamma
approximation of the law of nonlinear functionals of Gaussian fields,
Poisson measures and Rademacher sequences. In this respect, the most
relevant references are the following. In \cite{NouPecptrf}, Stein's
method and Malliavin calculus have been combined for the first time, in
the framework of the one-dimensional normal and Gamma approximations on
Wiener space. The findings of \cite{NouPecptrf} are extended in \cite
{NouPeAOP2} and \cite{NouPeRev}, dealing respectively with lower
bounds and multidimensional normal approximations. Reference~\cite
{NouPeReinPOIN} contains applications of the results of \cite
{NouPecptrf} to the derivation of second-order Poincar\'{e}
inequalities. References \cite{PSTU} and \cite{NouPeReinRAD} use
appropriate versions of the non-Gaussian Malliavin calculus in order to
deal with the one-dimensional normal approximation, respectively, of
functionals of Poisson measures and of functionals of infinite
Rademacher sequences.
Note that all the previously quoted references deal with the normal
and Gamma approximation of functionals of Gaussian fields, Poisson
measure and Rademacher sequences. The theory developed in the present
paper represents the first extension of the above quoted criteria to a
possibly non-Gaussian, non-Poisson and non-Rademacher framework.
\end{longlist}
\end{rem}

%s2 ###
\section{Wiener chaos}\label{S:WienerC}
In this section we briefly introduce the notion of (Gaussian) \textit
{Wiener chaos}, and point out some of its
crucial properties. The reader is referred to \cite{Nbook}, Chapter 1, or
\cite{Janson}, Chapter 2, for any unexplained
definition or result.
Let $\mathbf{G}=\{G_i \dvtx  i\geq1\}$ be a sequence of i.i.d. centered
Gaussian random variables with unit variance.

\begin{defi}\label{DefHER+Chaos}
1. The Hermite polynomials $\{H_q \dvtx  q\geq0\}$
are defined as $H_q = \delta^q\mathbf{1}$, where $\mathbf{1}$ is the function
constantly equal to 1, and $\delta$ is the divergence operator, acting
on smooth functions as $\delta f(x) = xf(x)-f'(x)$. For instance,
$H_0=1$, $H_1(x)=x$, $H_2(x) = x^2-1$, and so on. Recall that the class
$\{(q!)^{-1/2}H_q\dvtx q\geq0\}$ is an orthonormal basis of
$L^2(\mathbb{R}, (2\pi)^{-1/2}e^{-x^2/2}\,dx).$\vspace*{-6pt}
\begin{longlist}[(2)]
\item[2.] A \textit{multi-index} $q = \{q_i \dvtx  i\geq1 \}$ is a
sequence of
nonnegative integers such that $q_i \neq0$ only for a finite number
of indices $i$. We also write $\Lambda$ to indicate the class of all
multi-indices, and use the notation $|q| = \sum_{i\geq1} q_i$, for
every $q\in\Lambda$.
\item[3.] For every $d\geq0$, the $d$th \textit{Wiener chaos} associated
with $\mathbf{G}$ is defined as follows: $C_0 = \mathbb{R}$, and, for
$d\geq1$,
$C_d$ is the $L^2(P)$-closed vector space generated by random variables
of the type
$
\Phi(q)=\prod^{\infty}_{i=1} H_{q_i}(G_i)$, $q\in\Lambda$ and $|q|
= d$.
% \end{equation}
\end{longlist}
\end{defi}

\begin{example}
(i) The first Wiener chaos $C_1$ is the Gaussian space generated
by $\mathbf{G}$, that is, $F\in C_1$ if and only if
$ F = \sum_{i=1}^{\infty} \lambda_i G_i$ for some sequence $\{
\lambda_i
\dvtx  i\geq1\}\in\ell^2$.\vspace*{-6pt}
\begin{longlist}[(ii)]
\item[(ii)] Fix $d,N\geq2$ and let $f \dvtx  [N]^d \rightarrow
\mathbb{R}$ be
symmetric and vanishing on diagonals. Then,
an element of $C_d$ is, for instance, the following $d$-homogeneous sum:
%
%
%e16 ###
\begin{eqnarray}\label{QDXI}
Q_d(\mathbf{G}) &=& d! \sum_{\{i_1,\ldots,i_d\}\subset[N]^d}
f(i_1,\ldots,i_d) G_{i_1}\cdots G_{i_d}
\nonumber
\\[-8pt]
\\[-8pt]
\nonumber
&=& \sum_{1\leq
i_1,\ldots,i_d\leq N} f(i_1,\ldots,i_d) G_{i_1}\cdots G_{i_d}.
\end{eqnarray}
\end{longlist}
\end{example}

It is easily seen that two random variables belonging to a Wiener chaos
of different orders are orthogonal in $L^2(P)$. Moreover, since linear
combinations of polynomials are dense in $L^2(P,\sigma(\mathbf{G}))$, one
has that $L^2(P,\sigma(\mathbf{G}))=\bigoplus_{d\geq0} C_d$, that
is, any
square integrable functional of $\mathbf{G}$ can be written as an infinite
sum, converging in $L^2$ and such that the $d$th summand is an element
of $C_d$ [the \textit{Wiener--It\^{o} chaotic decomposition} of
$L^2(P,\sigma(\mathbf{G}))$]. It is often useful to encode the properties
of random variables in the spaces $C_d$ by using increasing
tensor powers of Hilbert spaces (see, e.g., \cite{Janson}, Appendix E,
for a collection of useful facts about tensor products).
To do this, introduce an (arbitrary) real separable Hilbert space
$\mathfrak{H}$ with scalar product $\langle\cdot, \cdot\rangle
_{\mathfrak{H}}$ and, for $d\geq2$, denote
by $\mathfrak{H}^{\otimes d}$ (resp. $\mathfrak{H}^{\odot d}$) the
$d$th tensor power (resp. symmetric
tensor power) of $\mathfrak{H}$; write, moreover, $\mathfrak
{H}^{\otimes0} =\mathfrak{H}^{\odot0}=
\mathbb{R}$ and $\mathfrak{H}^{\otimes1} =\mathfrak{H}^{\odot1}=
\mathfrak{H}$. Let $\{e_j \dvtx  j\geq1\}$ be an orthonormal basis of
$\mathfrak{H}$. With every multi-index
$q\in\Lambda$,\vspace*{-1pt} we associate the tensor $e(q)\in\mathfrak
{H}^{\otimes
|q|}$ given by
$
e(q) = e_{i_1}^{\otimes q_{i_1}}\otimes\cdots\otimes
e_{i_k}^{\otimes q_{i_k}},
$
where $\{q_{i_1},\ldots,q_{i_k}\}$ are the nonzero elements of~$q$. We
also denote by $\tilde{e}(q)\in
\mathfrak{H}^{\odot|q|}$ the canonical symmetrization of $e(q)$. It is
well known that, for every $d\geq2$,
the collection $\{\tilde{e}(q)\dvtx q\in\Lambda, |q|=d\}$ defines a
complete orthogonal system in
$\mathfrak{H}^{\odot d}$. For every $d\geq1$ and every $h\in
\mathfrak
{H}^{\odot d}$ with the form
$h = \sum_{q\in\Lambda, |q|=d}c_{q}\tilde{e}(q),$ we define
$I_d(h) = \sum_{q\in\Lambda, |q|=d}c_{q}\Phi(q).$
%where $\Phi(q)$ is given in (\ref{Eq:Q-chaos}).
We also recall that, for every $d\geq1$,
the mapping $I_d \dvtx  \mathfrak{H}^{\odot d}\rightarrow C_d $ %(as
%defined in (\ref{Eq:ISO}))
is onto, and provides an isomorphism between $C_d$ and the Hilbert
space $\mathfrak{H}^{\odot d}$, endowed with the norm $\sqrt{d!}\|
\cdot
\|_{\mathfrak{H}^{\otimes d}}$. In particular, for every $h,h'\in
\mathfrak{H}^{\odot d}$,
$
E[I_d(h)I_d(h')] = d! \langle h , h' \rangle_{\mathfrak{H}^{\otimes d}}.
$
If $\mathfrak{H} = L^2(A,\mathscr{A},\mu)$, with $\mu$ $\sigma$-finite
and nonatomic, then the operators $I_d$ are indeed (multiple)
Wiener--It\^{o} integrals.
\begin{example}\label{Ex:partsums}
By definition, $G_i = I_1(e_i)$, for every $i\geq1$.
Moreover, the random variable $Q_d(\mathbf{G})$ defined in (\ref{QDXI}) is
such that
%
%
%e17 ###
\begin{eqnarray}\label{EqWiener-Gausssum}
Q_d(\mathbf{G}) = I_d (h)\hspace*{180pt}
\nonumber
\\[-8pt]
\\[-8pt]
\eqntext{\mbox{where } h=d!
\displaystyle\sum
_{\{i_1,\ldots,i_d\}\subset[N]^d} f(i_1,\ldots,i_d) e_{i_1}\otimes\cdots\otimes e_{i_d}\in\mathfrak{H}^{\odot
d}.}
\end{eqnarray}
\end{example}

The notion of ``contraction'' is the key to prove the general bounds
stated in the forthcoming Section \ref{S:newbounds}.

\begin{defi}[(Contractions)]\label{Def:contraction}
Let $\{e_i \dvtx i \geq1 \}$ be a complete orthonormal system in $\HH
$, so
that, for every $m\geq2$, $\{e_{j_1}\otimes\cdots
\otimes
e_{j_m}\dvtx  j_1,\ldots,j_m\geq1\}$ is a complete orthonormal system in
$\HH^m$.
Let $f= \sum_{j_1,\ldots,j_p}a(j_1,\ldots,j_p)e_{j_1}\otimes\cdots
\otimes e_{j_p}\in\HH^{\odot p}$ and $g = \sum
_{k_1,\ldots,k_q}b(k_1,\ldots,k_q)e_{k_1}\otimes\cdots\otimes
e_{k_q}\in\HH^{\odot q}$, with\break $\sum_{j_1,\ldots,j_p}a(j_1,\ldots,j_p)^2
<\infty$ and $g = \sum_{k_1,\ldots,k_q}b(k_1,\ldots,k_q)^2 <\infty$ (note
that $a$ and $b$ need not vanish on diagonals). For every
$r=0,\ldots,p\wedge q$, the $r$th \textit{contraction} of $f$ and $g$ is
the element of $\HH^{\otimes(p+q-2r)}$ defined as
\begin{eqnarray*}
f\otimes_r g & =& \sum_{j_1,\ldots,j_{p-r}=1}^\infty\sum
_{k_1,\ldots,k_{q-r}=1}^{\infty} a\star_r b
(j_1,\ldots,j_{p-r},k_1,\ldots,k_{q-r})\\
&&\hspace*{95pt}{}\times e_{j_1} \otimes\cdots
\otimes e_{j_{p-r}}\otimes e_{k_1} \otimes\cdots
\otimes
e_{k_{q-r}} \\
&=& \sum_{i_1,\ldots,i_r=1}^\infty\langle f , e_{i_1}\otimes\cdot\cdot
\cdot\otimes e_{i_r}\rangle_{\HH^{\otimes r}}\otimes\langle g ,
e_{i_1}\otimes\cdots\otimes e_{i_r}\rangle_{\HH^{\otimes r}},
\end{eqnarray*}
where the kernel $a\star_r b$ is defined according to Definition \textup{\ref
{def-contr}}, by taking $N=\infty$.
\end{defi}

Plainly, $f\otimes_0 g=f\otimes g$ equals the tensor product of
$f$ and $g$ while, for $p=q$, $f\otimes_p g=\langle
f,g\rangle_{\HH^{\otimes p}}$. Note that, in general (and except
for trivial cases), the contraction $f\otimes_r g$ is \textit{not}
a symmetric element of $\HH^{\otimes(p+q-2r)}$. The canonical
symmetrization of $f\otimes_r g$ is written
$f\,\widetilde{\otimes}_r\, g$.
Contractions appear in multiplication formulae like the following one:

\begin{prop}[(Multiplication formulae)]\label{P:MultiplictForm}
If $f\in\HH^{\odot p}$ and
$g\in\HH^{\odot q}$, then
$
I_p(f)I_q(g)=\sum_{r=0}^{p\wedge q} r!\left({{p}\atop{r}}\right)\left({{q}\atop{r}}\right)
I_{p+q-2r}(f\,\widetilde{\otimes}_r\, g).
$
\end{prop}

Note that the previous statement implies that multiple integrals admit
finite moments of every order.
The next result (see \cite{Janson}, Theorem 5.10) establishes a more
precise property, namely, that random variables living in a finite sum
of Wiener chaos are hypercontractive.

\begin{prop}[(Hypercontractivity)]\label{P:Hyper} Let $d\geq1$\vspace*{1pt}
be a
finite integer and assume that
$F\in\bigoplus^d_{k=0} C_k$. Fix reals $2\leq p \leq q
<\infty$. Then
$
E[|F|^q]^{1/q} \leq(q-1)^{d/2} E[|F|^p]^{1/p}.
$
\end{prop}

%We will use hypercontractivity in order to deduce some of the bounds
%stated in the next section.

%s3 ###
\section{Normal and chi-square approximation on Wiener chaos} \label
{S:newbounds}
Starting from this section, and for the rest of the paper, we adopt the
following notation for distances between laws of real-valued random
variables. The symbol $d_{\mathrm{TV}}(F,G)$ indicates the \textit{total variation
distance} between the law of $F$ and $G$, obtained from~(\ref
{EQ:distance}) by taking $\mathscr{H}$ equal to the class of all indicators
of the Borel subsets of $\mathbb{R}$. The symbol $d_{\mathrm{W}}(F,G)$ denotes
the \textit{Wasserstein distance}, obtained from~(\ref{EQ:distance}) by choosing
$\mathscr{H}$ as
the class of all Lipschitz functions with Lipschitz constant less than
or equal to 1. The symbol $d_{\mathrm{BW}}(F,G)$ stands for the \textit{bounded
Wasserstein distance} (or \textit{Fortet--Mourier distance}),
deduced from (\ref{EQ:distance}) by choosing $\mathscr{H}$ as the
class of all Lipschitz functions that are bounded by 1, and
with Lipschitz constant less than or equal to 1. While $d_{\mathrm{Kol}}(F,G)
\leq d_{\mathrm{TV}}(F,G)$ and $d_{\mathrm{BW}}(F,G) \leq d_{\mathrm{W}}(F,G)$, in general,
$d_{\mathrm{TV}}(F,G)$ and $d_{\mathrm{W}}(F,G)$ are not comparable.

In what follows, we consider as given an i.i.d. centered standard
Gaussian sequence $ \mathbf{G} =\{G_i\dvtx  i \geq1\}$, and we shall
adopt the
Wiener chaos notation introduced in Section \ref{S:WienerC}.
%s3.1 ###
\subsection{Central limit theorems}
In the recent series of papers \cite
{NouPecptrf,NouPeAOP2,NouPeReinPOIN}, it has been shown that one
can effectively combine Malliavin calculus with Stein's method, in
order to evaluate the distance between
the law of an element of a fixed Wiener chaos, say, $F$, and a standard
Gaussian distribution. In this section
we state several refinements of these results, by showing, in
particular, that all the relevant bounds can be
expressed in terms of the fourth moment of $F$. The proof of the
following theorem involves the use of Malliavin calculus and is deferred
to Section~\ref{SS:Proof4thCUM}.

\begin{theorem}[(Fourth moment bounds)]\label{T:4thcumulant}
Fix $d\geq2$. Let $F=I_d(h)$, $h\in\HH^{\odot d}$, be an
element of
the $d$th Gaussian Wiener chaos $C_d$
such that $E(F^2)=1$,
let $Z\sim\mathscr{N}(0,1)$, and write
\begin{eqnarray*}
T_1(F) &:=& \sqrt{d^{2}\sum_{r=1}^{d-1}( r-1)
!^{2}\pmatrix{{d-1}\cr{r-1}}
^{4}( 2d-2r) !\Vert h\,\widetilde{\otimes
}_{r}\,h
\Vert_{
\mathfrak{H}^{\otimes2( d-r) }}^{2}}, \\
T_2(F) &:=& \sqrt{\frac{d-1}{3d}|E(F^4) - 3|}.
\end{eqnarray*}
We have
$
T_1(F)\leq T_2(F).
$
Moreover, $d_{\mathrm{TV}}(F,Z) \leq2T_1(F)$ and
$d_{\mathrm{W}}(F,Z)\leq T_1(F)$.
Finally, let $\varphi\dvtx  \mathbb{R}\rightarrow\mathbb{R}$ be a thrice
differentiable
function such that
$\|\varphi'''\|_{\infty} <\infty$. Then, one has that
$
|E[\varphi(F)] - E[\varphi(Z)] | \leq C_{*} \times T_1(F),
$
with
%
%
%e18 ###
\begin{eqnarray}\label{Cost2}
C_{*}&=& 4\sqrt{2}(1+5^{{3d}/2}
)
\nonumber
\\[-8pt]
\\[-8pt]
\nonumber
&&{}\times\max\biggl\{ \frac32\vert\varphi^{\prime\prime
}( 0)
\vert+\frac{\Vert\varphi^{\prime\prime\prime}
\Vert
_{\infty}}{3}\frac{2\sqrt{2}}{\sqrt{\pi}};2\vert\varphi
^{\prime
}( 0) \vert+\frac{1}{3}\Vert\varphi
^{\prime
\prime
\prime}\Vert_{\infty}\biggr\} .
\end{eqnarray}
\end{theorem}

\begin{rem}If $E(F)=0$
and $F$ has a finite fourth moment,
then the quantity $\kappa_4 (F)= E(F^4)-3E(F^2)^2$ is known as the
\textit{fourth cumulant} of $F$.
One can also prove (see, e.g., \cite{NP}) that, if $F$ is a nonzero
element of the $d$th Wiener chaos of a
given Gaussian sequence ($d\geq2$), then $\kappa_4(F) >0$.
\end{rem}

Now fix $d\geq2$, and consider a sequence of random variables of the
type $F_n = I_d(h_n)$, $n\geq1$,
such that, as $n\rightarrow\infty$, $E(F_n^2)=d!\|h_n\|^2_{\HH
^{\otimes
d}}\rightarrow1$. In \cite{NP}
it is proved that the following double implication holds: as
$n\rightarrow\infty$,
%
%
%e20 ###
%e19 ###
\begin{eqnarray}
&&\Vert h_n\,\widetilde{\otimes}_{r}\,h_n\Vert_{
\mathfrak{H}^{\otimes( 2d-2r) }} \rightarrow0\qquad
\forall
r=1,\ldots,d-1
\nonumber
\\[-8pt]
\\[-8pt]
\nonumber
&&\qquad\Leftrightarrow\quad  \Vert h_n \otimes_{r}
h_n\Vert_{
\mathfrak{H}^{\otimes( 2d-2r) }} \rightarrow0\qquad
\forall
r=1,\ldots,d-1.
\label{symcontr}
\end{eqnarray}
Theorem \ref{T:4thcumulant}, combined with (\ref{symcontr}), allows
therefore to recover the following characterization of CLTs on Wiener
chaos. It has been first proved (by
other methods) in \cite{NP}.

\begin{theorem}[(See \cite{NO,NP})]\label{T:NunuGio} Fix $d\geq2$, and
let $F_n = I_d (h_n)$, $n\geq1$ be a
sequence in the $d$th Wiener chaos of $\mathbf{G}$. Assume that $\lim
_{n\rightarrow\infty}E(F_n^2)= 1$. Then,
the following three conditions \textup{(1)}--\textup{(3)} are equivalent, as
$n\rightarrow\infty$:
\textup{(1)} $F_n$ converges in law to $Z\sim\mathscr{N}(0,1)$;
\textup{(2)} $E(F_n^4) \rightarrow E(Z^4)=3$;
\textup{(3)} for every $r=1,\ldots,d-1$, $\| h_n \otimes_{r} h_n\|
_{\mathfrak{H}^{\otimes2( d-r) }} \rightarrow0$.
\end{theorem}

\begin{pf}
Since $\sup_nE(F_n^2)<\infty$, one deduces from Proposition \ref
{P:Hyper} that, for every $M>2$, one has
$\sup_n E|F_n|^M<\infty$. By uniform integrability, it follows that, if
\textup{(1)} is in order, then necessarily
$E(F_n^4) \rightarrow E(Z^4)=3$. The rest of the proof is a consequence
of the bounds
in Theorem \ref{T:4thcumulant}.
\end{pf}

The following (elementary) result is one of the staples of the present paper.
%It relates the Hilbert space contractions of formula (\ref{defContr})
%with the discrete contraction operators defined in (
We state it in a form which is also useful for the chi-square
approximation of Section~\ref{SS:chisq}.

\begin{lemma}\label{Lem:ContNorms}Fix $d\geq2$, and suppose that
$h\in\HH^{\odot d}$ is given by
\textup{(\ref{EqWiener-Gausssum})}, with $f \dvtx  [N]^d \rightarrow\mathbb{R}$ symmetric
and vanishing on diagonals. Then, for $r=1,\ldots,d-1$,
$
\|h\otimes_r h\|_{\HH^{\otimes(2d-2r)}} = \|f\star_r f\|_{2d-2r},
$
where we have used the notation introduced in Definition \textup{\ref{def-contr}}.
Also, if $d$ is even, then, for every $\alpha_1,\alpha_2\in\mathbb{R}$,
$
\|\alpha_1 (h\otimes_{d/2} h) + \alpha_2 h \|_{\HH^{\otimes d}} = \|
\alpha_1 (f\star_{d/2} f )+ \alpha_2 f\|_{d}.
$
\end{lemma}
\begin{pf}Fix $r=1,\ldots,d-1$. Using (\ref{EqWiener-Gausssum}) and
the fact that $\{e_{j}\dvtx j\geq1 \}$ is an orthonormal basis of $\HH$,
one infers that
%
%
%e24 ###
%e23 ###
%e22 ###
%e21 ###
\begin{eqnarray}\label{omaracchia}
 h\otimes_r h &= &\sum_{1\leq i_1,\ldots,i_d\leq N }\sum
_{1\leq
j_1,\ldots,j_d \leq N} f(i_1,\ldots,i_d) f(j_1,\ldots,j_d) \nonumber\\
&&\hspace*{100pt}{}\times[e_{i_1}\otimes
\cdots\otimes e_{i_d}]\otimes_r [e_{j_1}\otimes\cdots\otimes
e_{j_d}] \nonumber\\
& =&\sum_{1\leq a_1,\ldots,a_r\leq N}\sum_{1\leq
k_1,\ldots,k_{2d-2r}\leq N}
f(a_1,\ldots,a_r,k_1,\ldots,k_{d-r})
\nonumber
\\[-8pt]
\\[-8pt]
\nonumber
&&{}\hspace*{120pt}\times f(a_1,\ldots,a_r,k_{d-r+1},\ldots,k_{2d-2r})
\\
&&{}\hspace*{120pt}\times e_{k_1}\otimes \cdots\otimes e_{k_{2d-2r}}\nonumber\\
&=& \sum_{1\leq k_1,\ldots,k_{2d-2r}\leq N}f\star_r
f(k_1,\ldots,k_{2d-2r})
e_{k_1}\otimes\cdots\otimes e_{k_{2d-2r}}.\nonumber
\end{eqnarray}
Since the set $\{e_{k_1}\otimes\cdots\otimes
e_{k_{2d-2r}}\dvtx k_1,\ldots,k_{2d-2r} \geq1 \}$ is an orthonormal basis
of $\HH^{\otimes(2d-2r)}$, one deduces immediately
$\|h\otimes_r h\|_{\HH^{\otimes(2d-2r)}} = \|f\star_r f\|_{2d-2r}$.
The proof of the other identity is analogous.
\end{pf}

\begin{rem} Theorem \ref{T:NunuGio} and Lemma \ref{Lem:ContNorms}
yield immediately a proof of Proposition \ref{T:NGCGHomSums} in the
case $m=1$.
\end{rem}

%s3.2 ###
\subsection{Chi-square limit theorems} \label{SS:chisq}

As demonstrated in \cite{NouPe08,NouPecptrf}, the combination of
Malliavin calculus and Stein's method also allows to estimate
the distance between the law of an element $F$ of a fixed Wiener chaos and
a (centered) chi-square distribution $\chi^2(\nu)$ with $\nu$ degrees
of freedom.
Analogously to the previous section for Gaussian approximations, we now
state a number of
refinements of the results proved in \cite{NouPe08,NouPecptrf}.
In particular, we will show that all the relevant bounds can be
expressed in terms of a specific linear combination of the third and fourth
moments of $F$. The proof is deferred to Section \ref{SS:Proof4thCUMbis}.

\begin{theorem}[(Third and fourth moment bounds)]\label{T:3and4}
Fix an even integer $d\geq2$ as well as an integer $\nu\geq1$.
Let $F=I_d(h)$ be an element of the $d$th Gaussian chaos $C_d$
such that $E(F^2)=2\nu$,
let $Z_\nu\sim\chi^2(\nu)$, and write
\begin{eqnarray*}
T_3(F) &:=& \biggl[
4d!\biggl\| h - \frac{d!^2}{4({d}/{2})!^3}h\,\widetilde{\otimes
}_{d/2}\,h\biggr\|^2_{\HH^{\otimes d}}\\
&&\hspace*{3pt}{}+d^{2}
\mathop{\sum_{r=1,\ldots,d-1}}_{r\neq d/2}
( r-1) !^{2} \pmatrix{{d-1}\cr{r-1}}
^{4} ( 2d-2r) !\Vert h\,\widetilde{\otimes
}_{r}\,h
\Vert_{
\mathfrak{H}^{\otimes2( d-r) }}^{2} \biggr]^{1/2}, \\
T_4(F) &:=& \sqrt{\frac{d-1}{3d}|E(F^4)-
12E(F^3)-12\nu
^2+48\nu|}.
\end{eqnarray*}
Then
$
T_3(F)\leq T_4(F)
$
and
$
d_{\mathrm{BW}}(F,Z_\nu) \leq
\max\{
\sqrt{\frac{2\pi}{\nu}},\frac1\nu+\frac2{\nu^2}
\}
T_3(F).
$
\end{theorem}

Now fix an even integer $d\geq2$, and consider a sequence of random
variables of the type
$F_n = I_d(h_n)$, $n\geq1$, such that, as $n\rightarrow\infty$,
$E(F_n^2)=d!\|h_n\|^2_{\HH^{\otimes d}}\rightarrow2\nu$.
In~\cite{NouPe08} it is proved that the following double implication holds:
as $n\rightarrow\infty$,
%
%
%e26 ###
%e25 ###
\begin{eqnarray}\label{symcontr2}
&& \Vert h_n\,\widetilde{\otimes}_{r}\,h_n\Vert_{
\mathfrak{H}^{\otimes2( d-r) }} \rightarrow0\qquad
\forall
r=1,\ldots,d-1, r\neq d/2
\nonumber
\\[-8pt]
\\[-8pt]
\nonumber
&&\qquad  \Longleftrightarrow\quad
\Vert h_n \otimes_{r} h_n\Vert_{
\mathfrak{H}^{\otimes2( d-r) }} \rightarrow0\qquad
\forall
r=1,\ldots,d-1, r\neq d/2.
\end{eqnarray}
Theorem \ref{T:3and4}, combined with (\ref{symcontr2}), allows
therefore to recover the following
characterization of chi-square limit theorems on Wiener chaos. Note
that this is a special case of a
``noncentral limit theorem''; one usually calls ``noncentral limit theorem''
any result involving convergence in law to a non-Gaussian distribution.

\begin{theorem}[(See \cite{NouPe08})]\label{T:IvanGio} Fix an even integer
$d\geq2$, and let $F_n = I_d (h_n)$, $n\geq1$ be a
sequence in the $d$th Wiener chaos of $\mathbf{G}$. Assume that $\lim
_{n\rightarrow\infty}E(F_n^2)= 2\nu$.
Then, the following three conditions \textup{(1)}--\textup{(3)} are
equivalent, as $n\rightarrow\infty$:
\textup{(1)} $F_n$ converges in law to $Z_\nu\sim\chi^2(\nu)$;
\textup{(2)} $E(F_n^4)-12E(F_n^3) \rightarrow E(Z_\nu^4)-12E(Z_\nu
^3)=12\nu
^2-48\nu$;
\textup{(3)} $\| h_n\,\widetilde{\otimes}_{d/2}\, h_n -
4(d/2)!^3d!^{-2}\times
h_n\|_{\HH^{\otimes d}}\to0$
and, for every $r=1,\ldots,d-1$ such that $r\neq d/2$,
$\| h_n \otimes_{r} h_n\|_{\mathfrak{H}^{\otimes2(
d-r)
}} \rightarrow0$.
\end{theorem}
\begin{pf}
Since $\sup_nE(F_n^2)<\infty$, one deduces from Proposition \ref
{P:Hyper} that, for every $M>2$, one has
$\sup_n E|F_n|^M<\infty$. By uniform integrability, it follows that, if
\textup{(1)} holds, then necessarily
$E(F_n^4)-12E(F_n^3) \rightarrow E(Z_\nu^4)-12E(Z_\nu^3)=12\nu
^2-48\nu$.
The rest of the proof is a
consequence of Theorem \ref{T:3and4}.
\end{pf}

\begin{rem}By using the second identity in Lemma \ref
{Lem:ContNorms} in the case $\alpha_1 =1$ and
$\alpha_2 = -4(\frac{d}2)!^3d!^{-2}$, Theorem \ref{T:IvanGio}
yields an immediate proof of Proposition \ref{T:IvGioNCLT}.
\end{rem}

%s4 ###
\section{Low influences and proximity of homogeneous sums}\label{S:MOO}

We now turn to some remarkable invariance principles by Rotar' \cite
{Rotar2} and Mossel, O'Donnell and Oleszkiewicz \cite{MOO}. As already
discussed, the results proved in \cite{Rotar2} yield sufficient
conditions in order to have that the laws of homogeneous sums (or, more
generally, polynomial forms) that are built from two different
sequences of independent random variables are asymptotically close,
whereas in \cite{MOO} one can find explicit upper bounds on the
distance between these laws. Since in this paper we adopt the
perspective of deducing general convergence results from limit theorems
on a Gaussian space, we will state the results of \cite{Rotar2} and
\cite{MOO} in a slightly less general form, namely, by assuming that
one of the sequences is i.i.d. Gaussian. See also Davydov and Rotar'
\cite{DavidovRotar}, and the references therein, for some general
characterizations of the asymptotic proximity of probability distributions.

\begin{theorem}[(See \cite{MOO})]\label{thmMOO}
Let $\mathbf{X}=\{X_i, i\geq1\}$ be a collection of
centered independent random variables with unit variance, and let
$\mathbf{G}=\{G_i\dvtx i\geq1\}$ be a collection of standard centered i.i.d.
Gaussian random variables.
Fix $d\geq1$, and let $\{N_n, f_n \dvtx  n\geq1\}$ be a sequence
such that
$\{N_n\dvtx  n\geq1\}$ is a sequence of integers
going to infinity, and each $f_n\dvtx  [N_n]^{d} \rightarrow\mathbb{R}$ is
symmetric and vanishes on diagonals.
Define $Q_{d}(N_n,f_n, \mathbf{X})$ and $Q_{d}(N_n,f_n, \mathbf{G})$
according to \textup{(\ref{EQ:HOMsums})}.
Recall the definition~\textup{(\ref{infi})} of $\operatorname{Inf}_i(f_n)$.
\begin{enumerate}
\item If $\sup_{i\geq1}E[|X_i|^{2+\varepsilon}]<\infty$ for some
$\varepsilon
>0$ and
if $\max_{1\leq i\leq N_n}
\operatorname{Inf}_i(f_n)
\to0$ as $n\to\infty$,
then
$\sup_{z\in\mathbb{R}}
|P[Q_{d}(N_n,f_n, \mathbf{X})\leq z]-P[Q_{d}(N_n,f_n,
\mathbf{G})\leq z]
|
\to0$
as $n\to\infty$.
\item If the random variables $X_i$ are identically distributed and
if
\[
\max_{1\leq i\leq N_n}
\operatorname{Inf}_i(f_n)
\to0\qquad \mbox{as }n\to\infty,
\]
then
$|E[\psi(Q_{d}(N_n,f_n, \mathbf{X}))]-E[\psi(Q_{d}(N_n,f_n,
\mathbf{G}))]|
\to0$ as $n\to\infty$,
for every continuous bounded function $\psi\dvtx  \mathbb{R}\rightarrow
\mathbb{R}$.
\item If $\beta:=\sup_{i\geq1}E[|X_i|^3]<\infty$, then, for
all thrice
differentiable $\varphi\dvtx \mathbb{R}\to\mathbb{R}$
such that $\|\varphi'''\|_\infty<\infty$ and for every fixed $n$,
$|E[\varphi(Q_{d}(N_n,f_n, \mathbf{X}))]-E[\varphi
(Q_{d}(N_n,f_n, \mathbf{G}))]|
\leq\|\varphi'''\|_\infty(30\beta)^d d! \sqrt{\max
_{1\leq i\leq N_n}
\operatorname{Inf}_i(f_n)}.$
\end{enumerate}
\end{theorem}
\begin{pf}
Point 1 is Theorem 2.2 in \cite{MOO}. Point 2 is Proposition 1 in
\cite
{Rotar2}. Point~3 is Theorem 3.18 (under Hypothesis H2) in \cite{MOO}.
Note that our polynomials $Q_d$ relate to polynomials $d! Q$ in \cite
{MOO}, hence the extra factor of $d!$ in
the bound.
\end{pf}

%associate explicit upper bounds (albeit non optimal) with the
%convergence (\ref{EQ:MOO1}).
%We stress that Point 3 in Theorem \ref{thmMOO} is a \textit{%non-asymptotic result}. In particular, the estimate
%in point 3 holds for every $n$ independently of the asymptotic
%behavior of the quantities $\max_{1\leq i\leq N_n}\operatorname{Inf}_i(f_n)$.
%and (
%}

In the sequel, we will also need the following technical lemma,
which follows directly by combining Propositions 3.11, 3.12 and 3.16 in
\cite{MOO}.
%, which is a combination
%of several results stated and proved in \cite{MOO}.
%
\begin{lemma}\label{hyper}
Let $\mathbf{X}=\{X_i, i\geq1\}$ be a collection of
centered independent random variables with unit variance. Assume,
moreover, that
$\gamma:=\sup_{i\geq1}E[|X_i|^q]<\infty$ for some $q>2$.
Fix $N,d\geq1$, and let $f\dvtx [N]^d\to\mathbb{R}$ be a symmetric function
(here, observe that we do not require that $f$ vanishes on
diagonals).
Define $Q_d(\mathbf{X})=Q_d(N,f,\mathbf{X})$ by \textup{(\ref{EQ:HOMsums})}.
Then $E[
|Q_d(\mathbf{X})|^q
]\leq
\gamma^d
(
2\sqrt{q-1})^{qd}\times
E[
Q_d(\mathbf{X})^2
]^{q/2}.$
\end{lemma}

As already evoked in the \hyperref[intr]{Introduction}, one of the key elements in the
proof of Theorem \ref{thmMOO} given in
\cite{MOO} is the use of an elegant probabilistic technique, which is
in turn inspired by the well-known
Lindeberg's proof of the central limit theorem.
%Very lucid accounts of Lindeberg's ideas can be found in Trotter
We will now state and prove a useful lemma, concerning moments of
homogeneous sums. We stress that the proof of
the forthcoming Lemma \ref{ivan1} could be directly deduced from the
general Lindeberg-type results
developed in \cite{MOO} (basically, by representing powers of
homogeneous sums as linear combinations of
homogeneous sums, and then by exploiting hypercontractivity). However,
this would require the introduction
of some more notation (in order to take into account different powers
of the same random variable), and we prefer to provide a direct proof,
which also serves as an illustration of some of the crucial techniques
of \cite{MOO}.\looseness=1
%As suggested before, we believe that the employed technique may replace
%``diagram formulae'' and similar combinatorial devices, that are
%customarily used in order to estimate moments of non-linear
%functionals of independent sequences
%(see e.g. \cite{PecTaqSURV,Sur}, as well as the estimates on moments
%of homogeneous sums given in de
%Jong's paper \cite{deJongMulti}).
%
\begin{lemma}\label{ivan1}
Let $\mathbf{X}=\{X_i\dvtx  i\geq1\}$ and $\mathbf{Y}=\{Y_i\dvtx  i\geq
1\}$
be two collections of centered independent random variables with unit variance.
Fix some integers~$N$, $d\geq1$, and let $f\dvtx [N]^d\to\mathbb{R}$ be
a symmetric
function vanishing on diagonals.
Define $Q_d(\mathbf{X})=Q_d(N,f,\mathbf{X})$ and $Q_d(\mathbf
{Y})=Q_d(N,f,\mathbf{Y})$
according to (\ref{EQ:HOMsums}).
\begin{enumerate}
\item Suppose $k\geq2$ is such that: \textup{(a)} $X_i$ and $Y_i$ belong to
$L^k(\Omega)$ for all $i\geq1$;
\textup{(b)}~$E(X_i^l)=E(Y_i^l)$ for all $i\geq1$ and $l\in\{2,\ldots,k\}$.
Then $Q_d(\mathbf{X})$ and $Q_d(\mathbf{Y})$ belong to $L^k(\Omega)$, and
$E[Q_d(\mathbf{X})^l]=E[Q_d(\mathbf{Y})^l]$
for all
$l\in\{2,\ldots,k\}$.
\item Suppose $m> k\geq2$ are such that: \textup{(a)} $\alpha:=\max
\{\sup
_{i\geq1}E|X_i|^m,\break
\sup_{i\geq1}E | Y_i|^m\}<\infty$; \textup{(b)}
$E(X_i^l)=E(Y_i^l)$ for all $i\geq1$ and $l\in\{2,\ldots,k\}$.
Assume, moreover, (for simplicity) that: \textup{(c)} $E[Q_d(\mathbf{X})^2]^{1/2}
\leq M$ for some finite constant \mbox{$M \geq1$}.
Then $Q_d(\mathbf{X})$ and $Q_d(\mathbf{Y})$ belong to $L^m(\Omega)$ and,
for all $l\in\{k+1,\ldots,m\}$,
$|E(Q_d(\mathbf{X})^l)-E(Q_d(\mathbf{Y})^l)|
\leq c_{d,l,m,\alpha} \times M^{l-k+1}
\times\max_{1\leq i\leq N}\{\max[
\operatorname{Inf}_i(f)
^{{k-1}/2};
\operatorname{Inf}_i(f)
^{{l}/2 -1}] \},
$
where\vspace*{1pt}
$c_{d,l,m,\alpha}=\break 2^{l+1}(d-1)!^{-1}\times\alpha^{{dl}/m}(2\sqrt
{l-1})^{(2d-1)l}
d!^{l-1}.$
\end{enumerate}
\end{lemma}

\begin{pf}
While Point 1 could be verified by a direct (elementary) computation,
we will obtain the same conclusion as
the by-product of a more sophisticated construction which will also
lead to the proof of Point 2. We shall assume, without loss of
generality, that the two sequences $\mathbf{X}$ and $\mathbf{Y}$ are
stochastically independent. For $i=0,\ldots,N$, let
$\mathbf{Z}^{(i)}$ denote the sequence $(Y_1,\ldots,Y_i,X_{i+1}, \ldots,X_N)$.
Fix a particular $i\in\{1,\ldots,N\}$, and write
\begin{eqnarray*}
U_i&=&\mathop{\sum_{1\leq i_1,\ldots,i_d\leq N}}_{\forall
k\dvtx  i_k\neq i}
f(i_1,\ldots,i_d)Z^{(i)}_{i_1}\cdots Z^{(i)}_{i_d},\\
V_i&=&\mathop{\sum_{1\leq i_1,\ldots,i_d\leq N}}_{\exists
k\dvtx  i_k= i}
f(i_1,\ldots,i_d)Z^{(i)}_{i_1}
\cdots\widehat{Z^{(i)}_i}\cdots Z^{(i)}_{i_d},
\end{eqnarray*}
where $\widehat{Z^{(i)}_i}$ means that this particular term is dropped
(observe that this notation bears no
ambiguity: indeed, since $f$ vanishes on diagonals, each string
$i_1,\ldots,i_d$ contributing to the definition of $V_i$ contains the
symbol $i$ exactly once).
Note that~$U_i$ and $V_i$ are independent of the variables $X_i$ and
$Y_i$, and that
$Q_d(\mathbf{Z}^{(i-1)})=U_i+X_iV_i$ and $Q_d(\mathbf{Z}^{(i)})=U_i+Y_iV_i$.
%By the binomial formula, we have that, for all $x,y\in\R$,
%(x+y)^l = \sum_{j=0}^{l} \binom{l}{j}x^{l-j}y^{j}.
By using the independence of $X_i$ and $Y_i$ from $U_i$ and $V_i$
[as well as the fact that $E(X_i^l)=E(Y_i^l)$ for all $i$ and all
$1\leq l\leq k$], we infer from the binomial formula that,
for $l\in\{2,\ldots,k\}$,
%
%
%e27 ###
\begin{eqnarray}\label{eqeq1}
E[ (U_i+X_iV_i)^l]
&=&\sum_{j=0}^{l} \pmatrix{{l}\cr{j}}E(U_i^{l-j}V_i^{j})E(X_i^{j})
\nonumber
\\[-8pt]
\\[-8pt]
\nonumber
&=&\sum_{j=0}^{l} \pmatrix{{l}\cr{j}}E(U_i^{l-j}V_i^{j})E(Y_i^{j})
=E[ (U_i+Y_iV_i)^l].
\end{eqnarray}
That is,
$E[ Q_d(\mathbf{Z}^{(i-1)})^l]=E[ Q_d(\mathbf{Z}^{(i)})^l]$
for all $i\in\{1,\ldots,N\}$ and $l\in\{2,\ldots,k\}$.
The desired conclusion of Point 1 follows by observing that $Q_d(\mathbf{
Z}^{(0)})= Q_d(\mathbf{X})
$ and $Q_d(\mathbf{Z}^{(N)})=Q_d(\mathbf{Y})$.
To prove Point 2, let $l\in\{k+1,\ldots,m\}$. Using (\ref{eqeq1}) and
then H\"older's inequality, we can write
\begin{eqnarray*}
&&\bigl|E\bigl[ Q_d\bigl(\mathbf{Z}^{(i-1)}\bigr)^l\bigr]-E\bigl[ Q_d\bigl(\mathbf{
Z}^{(i)}\bigr)^l
\bigr]\bigr|\\
&&\qquad=\Biggl|\sum_{j=k+1}^l \pmatrix{{l}\cr{j}}E(U_i^{l-j}V_i^{j})
\bigl(E(X_i^{j})-E(Y_i^j)\bigr)\Biggr|\\
&&\qquad\leq\sum_{j=k+1}^l \pmatrix{{l}\cr{j}}(E|U_i|^l)^{1-j/l}(E|V_i|^{l})^{j/l}
(
E|X_i|^{j}+E|Y_i|^j).
\end{eqnarray*}
By Lemma \ref{hyper}, since $E(U_i^2)\leq E(Q_d(\mathbf{
X})^2)\leq M ^2$,
we have
$
E|U_i|^l \leq\alpha^{dl/m}\times\break(2\sqrt{l-1}
)^{ld}E(U_i^2)^{l/2}
\leq\alpha^{dl/m}(2\sqrt{l-1})^{ld}M^l.
$
Similarly, since $E(V_i^2)=d!^2\times\break\operatorname{Inf}_i(f)$ [see (\ref{infi})], we have
$
E|V_i|^l \leq\alpha^{(d-1)l/m}(2\sqrt{l-1}
)^{l(d-1)}E(V_i^2)^{l/2}
\leq\break
\alpha^{(d-1)l/m}(2\sqrt{l-1})^{l(d-1)}d!^{l}(\operatorname{Inf}_i(f))^{l/2}.
$
Hence, since $E|Y_i|^j+E|X_i|^j\leq2\alpha^{j/m}$, we can write
\begin{eqnarray*}
&&\bigl|E\bigl[ Q_d\bigl(\mathbf{Z}^{(i-1)}\bigr)^l\bigr]-E\bigl[ Q_d\bigl(\mathbf{
Z}^{(i)}\bigr)^l\bigr]\bigr|\\
&&\qquad\leq2\sum_{j=k+1}^l\pmatrix{{l}\cr{j}}
\bigl(\alpha^{dl/m}\bigl(2\sqrt{l-1}\bigr)^{ld}M^l\bigr)^{1-{j}/l}\\
&&{}\hspace*{33pt}\qquad\quad \times
\bigl(
\alpha^{(d-1)/m}\bigl(2\sqrt{l-1}\bigr)^{d-1}d!\sqrt{\operatorname{Inf}_i(f)}
\bigr)^{j}
\alpha^{j/m}\\
&&\qquad\leq2^{l+1}\alpha^{{dl}/m}\bigl(2\sqrt{l-1}\bigr)^{l(2d-1)}
d!^{l}
M^{l-k-1}
\times\max\bigl[
\operatorname{Inf}_i(f)
^{{(k+1)}/2};
\operatorname{Inf}_i(f)
^{{l}/2}\bigr].
\end{eqnarray*}
Finally, summing for $i$ over $1,\ldots,N$ and using that\vspace*{-2pt} $\sum_{i=1}^N
\operatorname{Inf}_i(f)=\frac{\|f\|^2_d}{(d-1)!}\leq\frac{M^2}{d!(d-1)!}$
yields
\begin{eqnarray*}
&&|E[ Q_d(\mathbf{X})^l]-E[ Q_d(\mathbf{Y})^l]|\\
&&\qquad\leq 2^{l+1}\alpha^{{dl}/m}\bigl(2\sqrt
{l-1}\bigr)^{l(2d-1)}
d!^{l}
M^{l-k-1} \\
&&\qquad\quad{} \times\max_{1\leq i\leq N}\bigl\{
\max\bigl[
\operatorname{Inf}_i(f)
^{{(k-1)}/2};
\operatorname{Inf}_i(f)
^{{l}/2 -1}\bigr] \bigr\} \sum_{i=1}^N
\operatorname{Inf}_i(f) \\
&&\qquad\leq c_{d,l,m,\alpha} \times M^{l-k+1} \times\max
_{1\leq i\leq N}\bigl\{\max\bigl[
\operatorname{Inf}_i(f)
^{{(k-1)}/2};
\operatorname{Inf}_i(f)
^{{l}/2 -1}\bigr] \bigr\}.
\end{eqnarray*}\upqed
\end{pf}

%s5 ###
\section{Normal approximation of homogeneous sums}\label{S:normalHOMSUMS}
The following statement provides an explicit upper bound on the normal
approximation of homogenous
sums, when the test function has a bounded third derivative.
\begin{theorem}\label{mickey}
Let $\mathbf{X}=\{X_i, i\geq1\}$ be a collection of
centered independent random variables with unit variance.
Assume, moreover, that $\beta:=\sup_{i}E(X_i^4)<\infty$ and let $\alpha:=
\max
\{3; \beta\}$.
Fix $N,d\geq1$, and let $f\dvtx  [N]^{d} \rightarrow\mathbb{R}$ be
symmetric and
vanishing on diagonals.
Define $Q_d(\mathbf{X})=Q_{d}(N,f, \mathbf{X})$ according to \textup{(\ref
{EQ:HOMsums})} and assume that
$E[Q_d(\mathbf{X})^2]=1$.
Let $\varphi\dvtx \mathbb{R}\to\mathbb{R}$ be a thrice differentiable function
such that $\|\varphi'''\|_\infty\leq B$.
Then, for $Z\sim\mathscr{N}(0,1)$, we have, with $C_*$ defined by
\textup{(\ref{Cost2})},
%
%
%e29 ###
%e28 ###
\begin{eqnarray}\label{yes!}
&&| E[ \varphi(Q_d(\mathbf{X}))] - E[\varphi(Z)
]
|\nonumber\\
&&\qquad\leq B(30\beta)^dd!\sqrt{\max_{1\leq i\leq N}\operatorname{Inf}_i(f)}
\nonumber
\\[-8pt]
\\[-8pt]
\nonumber
&&\qquad\quad {}
+C_*\sqrt{\frac{d-1}{3d}}
\Bigl[\sqrt{|E[ Q_d(\mathbf{X})^4] -3|}
\\
&&\qquad{}\hspace*{74pt}+4\sqrt{2}\times
144^{d-1/2} \alpha^{{d}/2}\sqrt{d}d!
\Bigl(\max_{1\leq i\leq N}\operatorname{Inf}_i(f)\Bigr)^{1/4}
\Bigr].\nonumber
\end{eqnarray}
\end{theorem}
\begin{pf}
Let $\mathbf{G}=(G_i)_{i\geq1}$ be a standard centered i.i.d. Gaussian
sequence. We have
$| E[ \varphi(Q_d(\mathbf{X}))] - E[\varphi(Z)
]
|\leq
\delta_1+\delta_2$,
with
$\delta_1=| E[ \varphi(Q_d(\mathbf{X}))] - E[
\varphi
(Q_d(\mathbf{G}))]|$
and $\delta_2=| E[ \varphi(Q_d(\mathbf{G}))] - E[
\varphi
(Z)]|$.
By Theorem \ref{thmMOO}, we have
$\delta_1\leq B(30\beta)^d d!\sqrt{\max_{1\leq
i\leq N}\operatorname{Inf}_i(f)}$.
Since $E[ Q_d(\mathbf{X})^2]=E[ Q_d(\mathbf{G})^2]=1$,
Theorem \ref{T:4thcumulant} yields
$\delta_2\leq C_*\sqrt{\frac{d-1}{3d}|E[ Q_d(\mathbf{
G})^4]
-3|}.$
By Lemma \ref{ivan1}, Point 2\vspace*{-1pt} (with $M=1$, $k=2$ and $l=m=4$) and since
$\operatorname{Inf}_i(f)\leq1$ for all $i$, we have
$
|E[ Q_d(\mathbf{X})^4]
-E[ Q_d(\mathbf{G})^4]|\leq
32\times144^{2d-1} \alpha^d dd!^2
\sqrt{\max_{1\leq i\leq N}\operatorname{Inf}_i(f)},
$
so that
$
\delta_2\leq C_*\sqrt{\frac{d-1}{3d}}\times
[\sqrt{|E[ Q_d(\mathbf{X})^4] -3|}
+
4\sqrt{2}\times144^{d-1/2} \alpha^{{d}/2}\sqrt{d}d!
(\max_{1\leq i\leq N}\operatorname{Inf}_i(f))^{1/4}
].
$
%The proof of the theorem is done.
%Finally, the desired conclusion follows by putting the bounds for $
\end{pf}

\begin{rem}
 As a corollary of Theorem \ref{mickey}, we immediately recover de
Jong's Theorem \ref{T:deJong}, under the additional
hypothesis that $\sup_{i}E(X_i^4)<\infty$.
\end{rem}

As a converse statement, we now prove a slightly stronger version of
Theorem~\ref{T:byebyedeJong} stated in
Section \ref{sec-dejong}; an additional condition on contractions [see
assumption (5)
in Theorem \ref{T:byebyedeJongbis} just below and Definition \ref
{def-contr}]
has been added with respect to Theorem \ref{T:byebyedeJong}, making
the criterion more
easily applicable in practice.\looseness=1
\begin{theorem}\label{T:byebyedeJongbis}
We let the notation of Theorem \textup{\ref{T:byebyedeJong}} prevail.
Then, as $n\to\infty$, the assertions \textup{(1)--(4)} therein are
equivalent, and are also equivalent to
\textup{(5)} for all $r=1,\ldots,d-1$, $\|f_n\star_r f_n\|_{2d-2r}\to0$.
\end{theorem}
\begin{pf} The equivalences $\mathrm{(1)}\Leftrightarrow\mathrm{(2)} \Leftrightarrow
\mathrm{(5)}$ are a mere reformulation of
Theorem \ref{T:NunuGio}, deduced by taking into account the first
identity in Lemma \ref{Lem:ContNorms}.
On the other hand,
it is trivial that each one of conditions
\textup{(3)} and \textup{(4)} implies \textup{(1)}. So, it remains to prove the
implication $\mathrm{(1),(2),(5)}\Rightarrow\mathrm{(3),(4)}$.
Fix $z\in\mathbb{R}$. We have
$
| P[ Q_d(n,\mathbf{X})\leq z] - P
[Z\leq z]|
\leq| P[ Q_d(n,\mathbf{X})\leq z] - P
[ Q_d(n,\mathbf{
G})\leq z]|
+
| P[ Q_d(n,\mathbf{G})\leq z] - P[
Z\leq z]|
=:
\delta^{(a)}_n(z)+
\delta^{(b)}_n(z).
$
By assumption \textup{(2)} and Theorem \ref{T:4thcumulant}, we have
$\sup_{z\in\mathbb{R}}\delta^{(b)}_n(z)\to0$.
By combining assumption \textup{(5)} (for $r=d-1$) with (\ref{EQ:cruxcontr}),
we get that $\max_{1\leq i\leq N_n}\operatorname{Inf}_i(f_n)\to0$
as $n\to\infty
$. Hence,
Theorem \ref{thmMOO} (Point 1) implies that $\sup_{z\in\mathbb
{R}}\delta
^{(a)}_n(z)\to0$, and
the proof of the implication \textup{(1)}, \textup{(2)}, \textup{(5)}
$\Rightarrow$ \textup{(3)} is complete.
To prove that \textup{(1)} $\Rightarrow$ \textup{(4)}, one uses the same line
of reasoning,
the only difference being that we need to use Point 2 of Theorem \ref
{thmMOO} (along with the characterization of weak convergence based on
continuous bounded functions) instead of Point 1.
\end{pf}

Our techniques allow to directly control the Wasserstein distance between
the law of a homogenous sum and the law of a standard Gaussian random
variable, as illustrated by the
following result.

\begin{proposition}\label{mickeywasser}
As in Theorem \textup{\ref{mickey}}, let $\mathbf{X}=\{X_i, i\geq1\}$ be a
collection of
centered independent random variables with unit variance.
Assume, moreover, that $\beta:=\sup_{i}E(X_i^4)<\infty$ and note $\alpha
:=\max
\{ 3; \beta\}$.
Fix $N,d\geq1$, and let $f\dvtx  [N]^{d} \rightarrow\mathbb{R}$ be
symmetric and
vanishing on diagonals.
Define $Q_d(\mathbf{X})=Q_{d}(N,f, \mathbf{X})$ according to \textup{(\ref{EQ:HOMsums})}
and assume that $E[Q_d(\mathbf{X})^2]=1$.
Put
$
B_1 = 2 (30\beta)^d d!\sqrt{ \max_{1\leq i\leq N}\operatorname{Inf}_i(f)}
$
and\vspace*{-2pt}
$
B_2 = 12\sqrt{2} (
1+5^{{3d}/2}
) \sqrt{\frac{d-1}{3d}}
\times[ \sqrt{|E[ Q_d(\mathbf{X})^4] -3|}
+ 4\sqrt{2}\times144^{d-1/2} \alpha^{{d}/2}\sqrt{d}d!
(
{\max_{1\leq i\leq N}\operatorname{Inf}_i(f)})^{1/4}
].
$
For $Z\sim\mathscr{N}(0,1)$, we then have
$
d_W ( Q_d(\mathbf{X}), Z) \leq4 (B_1 + B_2)^{1/3},
$
provided $B_1+B_2 \leq\frac{3}{4 \sqrt{2}}$.
\end{proposition}

\begin{pf}
Let $h \in\operatorname{Lip}(1)$
be a Lipschitz function with constant 1.
By Radema\-cher's theorem, $h$ is Lebesgue-almost everywhere differentiable;
if we denote by $h'$ its derivative, then $\| h' \|_\infty\leq1$.
For $t > 0$, define
$h_t(x) = \int_{-\infty}^{\infty} h( \sqrt{t} y + \sqrt{1-t} x)
\phi
(y) \,dy,$
where $\phi$ denotes the standard normal density.
The triangle inequality gives
\begin{eqnarray*}
&& |E [h(Q_d(\mathbf{X}))] - E [ h(Z)]|\\
&&\qquad\leq | E [h_t(Q_d(\mathbf{X}))] - E [ h_t(Z)]| +
|
E [h(Q_d(\mathbf{X}))] - E [h_t(Q_d(\mathbf{X}))]|\\
&&\qquad\quad {}+ |E [h(Z)] - E [h_t(Z)]|.
\end{eqnarray*}
%
%Now, let us differentiate and integrate by parts, using that $\phi'(x)
%= - x \phi(x)$, to
%get
As
$
h_t''(x) = \frac{1-t}{\sqrt{t}}
\int_{-\infty}^{\infty} y h'( \sqrt{t} y + \sqrt{1-t} x) \phi
( y
) \,dy,
$ for $0 < t < 1$, we may bound
$
\| h_t'' \|_\infty\leq\frac{1-t}{\sqrt{t}} \| h' \|_\infty
\int
_{-\infty}^{\infty} \vert y \vert\phi(y) \,dy \leq\frac
{1}{\sqrt{t}}.
$
For $0 < t
\leq
\frac{1}{2}$ (so that $\sqrt{t} \leq\sqrt{1-t}$), we have
\begin{eqnarray*}
&&\vert E[h(
Q_d(\mathbf{X})
)] - E[h_t(Q_d(\mathbf{X}))]\vert\\
&&\qquad\leq \biggl\vert E\biggl[
\int_{-\infty}^{\infty} \bigl\{ h\bigl( \sqrt{t} y + \sqrt{1-t}
Q_d(\mathbf{
X})\bigr) - h\bigl(\sqrt{1-t} Q_d(\mathbf{X})\bigr)
\bigr\} \phi(y)\, dy\biggr] \biggr\vert\\
&&\qquad\quad {} + E\bigl[
\bigl\vert h\bigl( \sqrt{1-t} Q_d(\mathbf{X})\bigr) - h( Q_d(\mathbf{X}))
\bigr\vert\bigr]
\\
&&\qquad \leq \| h' \|_\infty\sqrt{t} \int_{-\infty}^{\infty}
| y | \phi
(y)\, dy +
\| h' \|_\infty\frac{t}{2 \sqrt{1-t}} E[ |Q_d(\mathbf{X})|]
\leq\frac32\sqrt{t}.
\end{eqnarray*}
Similarly,
$
\vert E [h(Z)] - E [h_t(Z)]\vert
\leq\frac{3}{2} \sqrt{t} .
$
We now apply Theorem \ref{mickey}. To bound $C_*$, we use that $\vert
h'_t(0) \vert\leq1$
and that $\vert h''_t(0) \vert\leq t^{-1/2}$;
also $\| h_t''' \|_\infty\leq2/t$ (as it can be shown by using the
same arguments as
above). Hence, as $2 \leq\frac{1}{t}$ and $\sqrt{2} \leq
t^{-1/2}$,
we have
\begin{eqnarray*}
C_*&\leq&4\sqrt{2}(1+5^{{3d}/{2}})
\times\max\biggl\{ \frac32
t^{-1/2} + \frac{4 \sqrt{2}}{3 \sqrt{\pi}} t^{-1}; 2 +
\frac{2}{3} t^{-1} \biggr\}\\
&\leq&4\sqrt{2}(1+5^{{3d}/{2}})\times\frac{3}{t}.
\end{eqnarray*}
Due to $\| h_t''' \|_\infty\leq2/t$, Theorem \ref{mickey}
gives the bound
$
| E[ h_t (Q_d(\mathbf{X}))] - E[h_t (Z)]|
\leq3 \sqrt{t} + (B_1 + B_2) \frac{1}{t}.
$
Minimizing $ 3 \sqrt{t} + (B_1 + B_2) \frac{1}{t} $ in $t$ gives that
$
t = ( \frac23 (B_1 + B_2) )^{2/3}$. Plugging in the values
and bounding the constant part ends the proof.\looseness=1
%\rightqed
\end{pf}

%s6 ###
\section{Chi-square approximation of homogeneous sums}\label{S:Chi2HOMSUMS}

The next result provides bounds on the chi-square approximation of
homogeneous sums.

\begin{theorem}\label{boundchi2}
Let $\mathbf{X}=\{X_i, i\geq1\}$ be a collection of
centered independent random variables with unit variance.
Assume, moreover, that
$\beta:=\sup_{i}E(X_i^4)<\infty$ and note $\alpha:=\max\{3;\beta\}$.
Fix an \textup{even} integer $d\geq2$ and, for $N\geq1$,
let $f\dvtx  [N]^{d} \rightarrow\mathbb{R}$ be symmetric and vanishing on
diagonals.
%For simplicity, assume that ${\rm Inf}_i(f)\leq1$ for all $1\leq i
Define $Q_d(\mathbf{X})=Q_{d}(N,f, \mathbf{X})$ according to \textup{(\ref{EQ:HOMsums})}
and assume that $E[Q_d(\mathbf{X})^2]=2\nu$ for some integer $\nu
\geq1$.
Let $\varphi\dvtx \mathbb{R}\to\mathbb{R}$ be a thrice differentiable function
such that $\|\varphi\|_\infty\leq1$, $\|\varphi'\|_\infty
\leq1$ and
$\|\varphi'''\|_\infty\leq B$.
Then, for $Z_\nu\sim\chi^2(\nu)$, we have
\begin{eqnarray*}
&&| E[ \varphi(Q_d(\mathbf{X}))] - E[\varphi(Z_\nu
)
]|\\
&&\qquad\leq B(30\beta)^dd!\sqrt{\max_{1\leq i\leq
N}\operatorname{Inf}_i(f)}\\
&&\qquad\quad {}+
\max\Biggl\{
\sqrt{\frac{2\pi}{\nu}},\frac1\nu+\frac2{\nu^2}
\Biggr\}\\
&&\qquad\qquad {}\times\Biggl(\sqrt{\frac{d-1}{3d}}
\Bigl[\sqrt{|E[ Q_d(\mathbf{X})^4] -12 E[ Q_d(\mathbf{
X})^3]-12\nu^2+48\nu|}\\
&&\hspace*{70pt}\qquad\quad {}+4\sqrt{d}d!\bigl(\sqrt2\times
144^{d-1/2} \alpha^{{d}/2}\\
&&\hspace*{150pt}{}+
\sqrt{\nu}\bigl(2\sqrt{2}\bigr)^{{3(2d-1)}/{2}}\alpha^{{3d}/2}\bigr)\\
&&\hspace*{175pt}\qquad\quad {}\times\Bigl(\max_{1\leq i\leq N}\operatorname{Inf}_i(f)\Bigr)^{1/4}
\Bigr]\Biggr).
\end{eqnarray*}
\end{theorem}

\begin{pf}
We proceed as in Theorem \ref{mickey}.
Let $\mathbf{G}=(G_i)_{i\geq1}$ denote a standard centered i.i.d.
Gaussian sequence. We have
$
| E[ \varphi(Q_d(\mathbf{X}))] - E[\varphi(Z_\nu
)
]|\leq
\delta_1+\delta_2
$
with
$
\delta_1=| E[ \varphi(Q_d(\mathbf{X}))] - E[
\varphi
(Q_d(\mathbf{G}))]|
$
and
$\delta_2=| E[ \varphi(Q_d(\mathbf{G}))] - E[
\varphi
(Z_\nu)]|.
$
By Theorem \ref{thmMOO} (Point 3), we have
$\delta_1\leq B(30\beta)^dd!\sqrt{\max_{1\leq
i\leq N}\operatorname{Inf}_i(f)}$.
By Theorem \ref{T:3and4}, we have, with $C_\#=\max\{
\sqrt{\frac{2\pi}{\nu}},\frac1\nu+\frac2{\nu^2}
\}$, that
$
(\delta_2)^2\leq
%d_{\mathrm{BW}}(Q_d(\mathbf{G}),Z_\nu)\leq
(C_\#)^2\times
\frac{d-1}{3d}|E[ Q_d(\mathbf{G})^4] - 12 E[Q_d(\mathbf{
G})^3]
-12\nu^2 + 48\nu|.
$
Additionally to the bound for $|E[ Q_d(\mathbf{X})^4]
-E[ Q_d(\mathbf{G})^4]|$ in Theorem \ref{mickey}, we have,
by Lemma \ref{ivan1},
$
|E[ Q_d(\mathbf{X})^3]
-E[ Q_d(\mathbf{G})^3]|
\leq
16\nu(2\sqrt{2})^{3(2d-1)}\alpha^{{3d}/4}dd!
\sqrt{\max_{1\leq i\leq
N}\operatorname{Inf}_i(f)}.
$\\
Hence, the proof is concluded since
\begin{eqnarray*}
&&\delta_2\leq C_\#\sqrt{\frac{d-1}{3d}}
\Bigl[\sqrt{|E[ Q_d(\mathbf{X})^4] -12 E[ Q_d(\mathbf{
X})^3]-12\nu^2+48\nu|}\\
&&\hspace*{51pt}\qquad{}+4\sqrt{d}d!\bigl(\sqrt2
\times144^{d-1/2} \alpha^{{d}/2}\\
&&\hspace*{96pt}\qquad{}+ \sqrt{\nu} \bigl(2\sqrt
{2}\bigr)^{{3(2d-1)}/{2}}
\alpha^{{3d}/2}\bigr)\\
&&\hspace*{161pt}\qquad{}\times\Bigl(\max_{1\leq i\leq N}\operatorname{Inf}_i(f)
\Bigr)^{1/4}\Bigr].
\end{eqnarray*}
\upqed\end{pf}

As an immediate corollary of Theorem \ref{boundchi2}, we deduce the
following new criterion for the asymptotic nonnormality of homogenous
sums---compare with Theorem \ref{T:deJong}.
\begin{cor}\label{T:deJongchi2} Let $\mathbf{X } = \{X_i \dvtx  i\geq
1 \}$
be a sequence of independent
centered random variables with unit variance such that $\sup_i E(X_i^4)
< \infty$.
Fix an \textup{even} integer $d\geq2$, and let $\{N_n, f_n \dvtx
n\geq1\}$
be a sequence such that
$\{N_n\dvtx  n\geq1\}$ is a sequence of integers
going to infinity, and each $f_n\dvtx  [N_n]^{d} \rightarrow\mathbb{R}$ is
symmetric and vanishes on diagonals.
Define $Q_d(n,\mathbf{X})=Q_{d}(N_n,f_n, \mathbf{X})$ according to \textup{(\ref
{EQ:HOMsums})}.
If, as $n\rightarrow\infty$,
$\textup{(i)}$ $E(Q_{d}(n, \mathbf{X})^2)\to2\nu$;
$\textup{(ii)}$ $E[Q_{d}(n, \mathbf{X})^4]- 12E[Q_{d}(N_n,f_n, \mathbf{
X})^3]\rightarrow12\nu^2-48\nu$; and
$\textup{(iii)}$ $\max_{1\leq i\leq N_n} \operatorname{Inf}_i(f_n)
\rightarrow0$;
then $Q_{d}(n, \mathbf{X})$ converges in law to $Z_\nu\sim\chi^2(\nu)$.
\end{cor}

The following statement contains a universal chi-square limit theorem
result: it is a general version of
Theorem \ref{T:seeyadeJong}.

\begin{theorem}\label{T:byebyedeJongter}
We let the notation of Theorem \textup{\ref{T:seeyadeJong}} prevail.
Then, as $n\to\infty$, the assertions \textup{(1)--(4)} therein are
equivalent, and are also equivalent to
\textup{(5)} $\| f_n\,\widetilde{\star}_{d/2}\, f_n - 4(d/2)!^3d!^{-2}\times
f_n\|_{d}\to0$
and, for every $r=1,\ldots,d-1$ such that $r\neq d/2$,
$\| f_n \star_{r} f_n\|_{2d-2r} \rightarrow0$.
\end{theorem}
\begin{pf}
The proof follows exactly the same lines of reasoning as in Theorem~\ref{T:byebyedeJongbis}. Details are left to the reader.
Let us just mention that the only differences consist in the use of
Theorem \ref{T:IvanGio} instead of
Theorem \ref{T:NunuGio}, and the use of Theorem~\ref{T:3and4}
instead of
Theorem \ref{T:4thcumulant}.
\end{pf}

%s7 ###
\section{Multivariate extensions}\label{S:multivariateext}

%s7.1 ###
\subsection{Bounds}
We recall here the standard multi-index notation. A multi-index is a
vector $\alpha\in\{0,1,\ldots\}^m$.\vspace*{-3pt}
We write
$
|\alpha|=\sum_{j=1}^m \alpha_j$,
$\alpha!=\prod_{j=1}^m \alpha_j!$,
$\partial_j =\frac{\partial}{\partial x_j}$,
$\partial^\alpha=\partial_1^{\alpha_1}\,\cdots\,\partial_d^{\alpha_d}$,
and $x^\alpha= \prod_{j=1}^m x_j^{\alpha_j}.
$\vspace*{1pt}
Note that, by convention, $0^0=1$. Also note that $|x^\alpha|=y^\alpha
$, where $y_j=|x_j|$ for all $j$.
Finally, for $\varphi\dvtx \mathbb{R}^m\to\mathbb{R}$ regular and
$k\geq1$, we put
$
\|\varphi^{(k)}\|_\infty=\max_{|\alpha|=k}\frac1{\alpha!}\sup
_{z\in\mathbb{R}
^m}|\partial^\alpha\varphi(z)|.
$

The forthcoming Theorem \ref{T:MULTIMOO} is a multivariate version of
Theorem \ref{thmMOO} (Point~3).
Observe that its statement (and its proof as well) follows closely
(\cite{Mossel}, Theorem 4.1). However, the result of \cite{Mossel}
is stated and proved under the assumption that one of the two i.i.d.
sequences lives
on a discrete probability space, hence, a~bit more work is needed.

%Since our framework also involves continuous i.i.d.
%sequences, we decided to provide a complete and self-contained proof.

%the proof is omitted, as it can be achieved along the lines of

\begin{theorem}\label{T:MULTIMOO}
Let $\mathbf{X}=\{X_i, i\geq1\}$ be a collection of
centered independent random variables with unit variance and such that
$\beta:=\sup_{i\geq1}E[|X_i|^3]<\infty$.
Let $\mathbf{G}=\{G_i \dvtx  i\geq1\}$ be a
standard centered i.i.d. Gaussian sequence.
Fix integers $m\geq1$, $d_m\geq\cdots\geq
d_1\geq1$ and $N_1,\ldots
,N_m\geq1$.
For every $j=1,\ldots,m$, let
$f_j\dvtx  [N_j]^{d_j} \rightarrow\mathbb{R}$ be a symmetric function
vanishing on
diagonals.
Define $Q^j(\mathbf{G})=Q_{d_j}(N_j,f_j, \mathbf{G})$ and $Q^j(\mathbf{
X})=Q_{d_j}(N_j,f_j, \mathbf{X})$
according to \textup{(\ref{EQ:HOMsums})}, and assume that $E[Q^j(\mathbf{
G})^2]=E[Q^j(\mathbf{X})^2]=1$ for all $j=1,\ldots,m$.
Assume that there exists a $C>0$ such that
$ \sum_{i=1}^{\max_j N_j } \max_{1\leq j \leq m} \operatorname{Inf}_i(f_j)\leq C$.
Then, for all thrice differentiable $\varphi\dvtx \mathbb{R}^m\to\mathbb{R}$
with $\|\varphi'''\|_\infty<\infty$,
we have
\begin{eqnarray*}
&&|E[\varphi(Q^1(\mathbf{X}),\ldots,Q^m(\mathbf{X}))]-E[
\varphi(Q^1(\mathbf{G}),\ldots,Q^m(\mathbf{G}))]|\\
&&\qquad\leq C\|\varphi'''\|_\infty\Biggl(\beta+ \sqrt{\frac8\pi
}\Biggr)
\Biggl[
\sum_{j=1}^m \bigl(16\sqrt2\beta\bigr)^{(d_j-1)/3}d_j!
\Biggr]^3 \\
&&\quad \qquad{}\times\sqrt{ \max_{1\leq j \leq m} \max
_{1\leq i\leq\max_j N_j}
\operatorname{Inf}_i(f_j)}.
\end{eqnarray*}
\end{theorem}

Observe that, in the one-dimensional case ($m=1$),
\[
 \sum_{i=1}^{\max_j
N_j } {\max_{1\leq j \leq m}} \operatorname{Inf}_i(f_j) = [d!(d-1)!]^{-1},
\]
so we can choose $C= [d!(d-1)!]^{-1}$. In this case, when $\beta$ is large,
the bound from
Theorem \ref{T:MULTIMOO} essentially differs from the one in Theorem
\ref{thmMOO} by a constant times a factor $d$.

\begin{pf*}{Proof of Theorem \protect\ref{T:MULTIMOO}}
Abbreviate
$
\mathbf{{Q}}(\mathbf{X}) = (Q^1(\mathbf{X}),\ldots,Q^m(\mathbf{X})),
$
and define $\mathbf{{Q}}(\mathbf{G})$ analogously. We proceed as for Lemma
\ref{ivan1}, with similar notation.
For $i=0,\ldots, \max_j N_j$, let
$\mathbf{Z}^{(i)}$ denote the sequence $(G_1,\ldots,G_i,X_{i+1},\break
\ldots
, X_{\max_j N_j } )$.
Using the triangle inequality,
\[
| E[ \varphi( \mathbf{{Q}}(\mathbf{X}) )] - E[ \varphi
( \mathbf{
{Q}}(\mathbf{G}) )]|
\leq \sum_{i=1}^{\max_j N_j }
\bigl|E\bigl[ \varphi\bigl( \mathbf{{Q}}\bigl(\mathbf{Z}^{(i-1)}\bigr)\bigr)\bigr] - E\bigl[
\varphi
\bigl(\mathbf{{Q}}\bigl(\mathbf{Z}^{(i)}\bigr)\bigr)\bigr]\bigr|.
\]
Now we can proceed as for inequality (31) in the proof of \cite{Mossel}, Theorem
4.1 to obtain
%Fix a particular $i\in\{1,\ldots, \max_j N_j \}$, and write, for $j=1,
%U_i^{(j)}= \sum_{\substack{1\leq i_1,\ldots,i_d\leq N_j\\
% f_j ( i_1,\ldots,i_d)Z^{(i)}_{i_1}\ldots Z^{(i)}_{i_d},
%V_i^{(j)}= \sum_{\substack{1\leq i_1,\ldots,i_d\leq N_j\\
% f_j (i_1,\ldots,i_d)Z^{(i)}_{i_1}
%%\ldots\widehat{Z^{(i)}_i}\ldots Z^{(i)}_{i_d},
%where $\widehat{Z^{(i)}_i}$ means that this particular term is dropped.
%We write
%$\mathbf{U}_i= ( U_i^{(1)}, \ldots, U_i^{(m)} )$ and $\mathbf{V}_i= (
%V_i^{(1)}, \ldots, V_i^{(m)} )$.
%Note that $\mathbf{U}_i$ and $\mathbf{V}_i$ are independent of the
%variables
%$X_i$ and $G_i$, and that
%Z}^{(i)})=\mathbf{U}_i+G_i \mathbf{V}_i$.
%y Taylor's theorem, using the independence of $X_i$ from $\mathbf{U}_i$
%and $\mathbf{V}_i$
%and $E(X_i)=0$ and $E(X_i^2)=1$, we have
%&&| E[ \varphi(\mathbf{U}_i+X_i\mathbf{V}_i)]-
%| \\
%&\leq&
% \sum_{|\alpha| = 3} \frac1{\alpha!} \sup_{z \in\R^m} |
%Similarly,
%$
%| E[ \varphi(\mathbf{U}_i+G_i\mathbf{V}_i)]-
% {\sum_{ \vert\alpha\vert\leq2 }}
%|
%$
%Due to the independence and the matching moments up to second order,
%we obtain that
%
\begin{eqnarray*}
&&\bigl|E\bigl[ \varphi\bigl(\mathbf{Q} \bigl(n,\mathbf{Z}^{(i-1)}\bigr)\bigr)\bigr] - E\bigl[
\varphi\bigl( \mathbf{Q} \bigl(n,\mathbf{Z}^{(i)}\bigr)\bigr)\bigr]\bigr|\\
&&\qquad=|E[ \varphi(\mathbf{U}_i+X_i \mathbf{V}_i)] - E[
\varphi(
\mathbf{U}_i+G_i \mathbf{V}_i)]|\\
&&\qquad\leq \Biggl(\beta+\sqrt{\frac8\pi}\Biggr)
\|\varphi'''\|_\infty
\sum_{|\alpha| = 3} E( |\mathbf{V}_i^\alpha|).
\end{eqnarray*}
While \cite{Mossel}, Theorem 4.1, now uses hypercontractivity results
for random variables on
finite probability spaces, here we bound the moments directly.
Abbreviate
$
\tau_{i} =\max_{1\leq j \leq m}
\operatorname{Inf}_i(f_j).
$
Next we use that, for $j=1, \ldots, m$, by Lemma \ref{hyper} (with
$q=3$), we have
$
E[ |V_i^{(j)}|^3]\leq(16\sqrt2\beta)^{d_j-1} E[ (V_i^{(j)}
)^2]^{3/2} = (16\sqrt2\beta)^{d_j-1}
d_j!^3 \tau_{i}^{3/2}.
$
Thus,
\begin{eqnarray*}
\sum_{ \vert\alpha\vert= 3}
E | (\mathbf{V}_i)^\alpha|
&=& \sum_{j,k,l=1}^m E \bigl(\bigl|V_i^{(j)} V_i^{(k)} V_i^{(l)}
\bigr|\bigr)\\
&\leq&\sum_{j,k,l=1}^m E \bigl(\bigl|V_i^{(j)}\bigr|^3\bigr)^{1/3} E
\bigl(\bigl|V_i^{(k)}\bigr|^3\bigr)^{1/3} E \bigl(\bigl|V_i^{(l)}\bigr|^3\bigr)^{1/3}\\
&=&\Biggl(\sum_{j=1}^m E
\bigl(\bigl|V_i^{(j)}\bigr|^3\bigr)^{1/3}\Biggr)^3\\
&\leq&\Biggl[\sum_{j=1}^m \bigl(16\sqrt2\beta\bigr)^{(d_j-1)/3}d_j!\Biggr]^3
\tau_{i}^{3/2}.
\end{eqnarray*}
Collecting the bounds, summing over $i$, and using that
$\sum_{i=1}^{\max_j N_j} \tau_{i} \leq C$ gives
the desired result.
\end{pf*}

The next statement gives explicit bounds on the distance to the normal
distribution for the
distribution of the vector
$(Q^1(\mathbf{X}),\ldots,Q^m(\mathbf{X}))$.

\begin{theorem} \label{multimoo2}
Let $\mathbf{X}=\{X_i \dvtx  i\geq1\}$ be a collection of centered independent
random variables
with unit variance. Assume, moreover, that $\beta: =\sup_i E
[|X_i|^3]<\infty$.
Fix integers $m\geq1$, $d_m\geq\cdots\geq
d_1\geq2$ and $N_1,\ldots
,N_m\geq1$.
For every $j=1,\ldots,m$, let
$f_j\dvtx  [N_j]^{d_j} \rightarrow\mathbb{R}$ be a symmetric function
vanishing on
diagonals.
Define $Q^j(\mathbf{X})=Q_{d_j}(N_j,f_j, \mathbf{X})$ according to\vspace*{-3pt} \textup{(\ref
{EQ:HOMsums})},
and assume that
$E[{Q^j(\mathbf{X})}^2]=1$ for all $j=1,\ldots,m$.
Let $V$ be the $m\times m$ symmetric matrix
given by
$V(i,j)=E[Q^i(\mathbf{X})Q^j(\mathbf{X})].$
Let $C$ be as in Theorem \textup{\ref{T:MULTIMOO}}.
Let $\varphi\dvtx \mathbb{R}^m\to\mathbb{R}$ be a thrice differentiable
function such that
$\|\varphi''\|_\infty<\infty$ and $\|\varphi'''\|_\infty<\infty$.
Then, for $Z_V=(Z^1_V,\ldots,Z^m_V)\sim\mathscr{N}_m(0,V)$
(centered Gaussian vector with covariance matrix $V$), we have
\begin{eqnarray*}
&&|E[\varphi(Q^1(\mathbf{X}),\ldots,Q^m(\mathbf{X}))]-E[
\varphi(Z_V)]|\\
&&\qquad\leq\|\varphi''\|_\infty
\Biggl(\sum_{i=1}^m\Delta_{ii}+2\sum_{1\leq i<j\leq m}
\Delta_{ij}\Biggr)\\
&&\quad \qquad{}+C\|\varphi'''\|_\infty\Biggl(\beta+ \sqrt{\frac8\pi}\Biggr)
\Biggl[
\sum_{j=1}^m \bigl(16\sqrt2\beta\bigr)^{(d_j-1)/3}d_j!
\Biggr]^3\\
&&\qquad\qquad{}\times \sqrt{\max_{1\leq j\leq m}\max_{1\leq
i\leq N_j}\operatorname{Inf}_i(f_j)}
\end{eqnarray*}
for $\Delta_{ij}$ given by
%
%
%e32 ###
%e31 ###
%e30 ###
\begin{eqnarray}\label{deltaij}
&&\frac{d_j}{\sqrt2}\sum_{r=1}^{d_i-1}(r-1)!\pmatrix{{d_i-1}\cr{r-1}}\pmatrix{
{d_j-1}\cr{r-1}}\nonumber\\
&&\qquad{}\quad\times\sqrt{(d_i+d_j-2r)!}
(
\|f_i\star_{d_i-r}f_i\|_{2r}
+\|f_j\star_{d_j-r}f_j\|_{2r}
)
\\
&&\qquad\quad {}+\mathbf{1}_{\{d_i<d_j\}}\sqrt{d_j!\pmatrix{{d_j}\cr{d_i}}\|f_j\star
_{d_j-d_i}f_j\|_{2d_i}}.\nonumber
\end{eqnarray}
\end{theorem}
\begin{pf}
The proof is divided into four steps.

\textit{Step 1: Reduction of the problem}. Let $\mathbf
{G}=(G_i)_{i\geq1}$
be a standard centered~i.i.d. Gaussian sequence. We have
$
|E[\varphi(Q^1(\mathbf{X}),\ldots,Q^m(\mathbf{X}))]-E[
\varphi(Z_V)]|
\leq
\delta_1+\delta_2
$
with
$
\delta_1=
| E[\varphi(Q^1(\mathbf{X}),\ldots,Q^m(\mathbf{X}))] -
E
[\varphi(Q^1(\mathbf{G}),\ldots,Q^m(\mathbf{G}))]|$
and $\delta_2=| E[\varphi(Q^1(\mathbf{G}),\ldots,Q^m(\mathbf{
G}))]-
E[ \varphi(Z_V)]|$.

\textit{Step 2: Bounding $\delta_1$}.
By Theorem \ref{T:MULTIMOO}, we have
\[
\delta_1\leq C\|\varphi'''\|_\infty\Biggl(\beta+ \sqrt{\frac
8\pi}\Biggr)
\Biggl[
\sum_{j=1}^m \bigl(16\sqrt2\beta\bigr)^{{(d_j-1)}/3}d_j!
\Biggr]^3 \sqrt{\max_{1\leq j\leq m}\max_{1\leq
i\leq N_j}\operatorname{Inf}_i(f_j)}
.
\]

\textit{Step 3: Bounding $\delta_2$}.
We will not use the result proved in \cite{NouPeRev},
since here we do not assume that the matrix $V$ is positive definite.
Instead, we will rather use an interpolation technique.\vspace*{-1pt} %{\it\`a la}
%Talagrand \cite{tala}.
Without loss of generality, we assume in this step
that $Z_V$ is independent of $\mathbf{G}$.
By (\ref{EqWiener-Gausssum}), we have that
$\{Q^j(\mathbf{G})\}_{1\leq j\leq m} \stackrel{\mathrm{Law}}{=} \{I_{d_j}
(h_j)\}_{1\leq j\leq m}$ where
$
h_j= {d_j! \sum_{\{i_1,\ldots,i_{d_j}\}\subset
[N_j]^{d_j}} } f_j(i_1,\ldots,i_{d_j}) e_{i_1}\otimes\cdots
\otimes e_{i_{d_j}}\in\mathfrak{H}^{\odot d},
$\vspace*{-3pt}
with $\HH=L^2([0,1])$ and $\{e_j\}_{j\geq1}$ any orthonormal
basis of
$\HH$.
For $t\in[0,1]$, set
$\Psi(t)=E[\varphi(\sqrt{1-t}(I_{d_1} (h_1),\ldots,I_{d_m}
(h_m))+\sqrt{t}Z_V)],$
so that
$\delta_2=|\Psi(1)
-\Psi(0)|\leq\sup_{t\in(0,1)}|\Psi'(t)|.$
We easily see that $\Psi'(t) =\break
\sum_{i=1}^m E [\frac{\partial\varphi}{\partial x_i}
(\sqrt
{1-t}(I_{d_1} (h_1),\ldots, I_{d_m} (h_m))
+\sqrt{t}Z_V)(
\frac1{2\sqrt{t}}Z^i_V-\frac{1}{2\sqrt{1-t}}I_{d_i} (h_i))
].
$
By integrating by parts, we can write
\begin{eqnarray*}
&&E\biggl[
\frac{\partial\varphi}{\partial x_i}\bigl(\sqrt{1-t}(I_{d_1}
(h_1),\ldots
,I_{d_m} (h_m))
+\sqrt{t}Z_V\bigr)Z_V^i
\biggr]\\
&&\qquad=\sqrt{t}\sum_{j=1}^mV(i,j) E\biggl[
\frac{\partial^2\varphi}{\partial x_i\,\partial x_j}\bigl(\sqrt
{1-t}(I_{d_1} (h_1),\ldots,I_{d_m} (h_m))
+\sqrt{t}Z_V\bigr)
\biggr].
\end{eqnarray*}
By using (\ref{Mexico2}) below in order to perform the integration by
parts, we can also write\looseness=1
\begin{eqnarray*}
&&E\biggl[
\frac{\partial\varphi}{\partial x_i}\bigl(\sqrt{1-t}(I_{d_1}
(h_1),\ldots
,I_{d_m} (h_m))
+\sqrt{t}Z_V\bigr)I_{d_i} (h_i)
\biggr]\\
%&=&E\{
%E[
%+\sqrt{t}z)I_{d_i} (h_i)
%]_{|z=Z_V}\}\\
%&=&\frac{\sqrt{1-t}}{d_i}\sum_{j=1}^m E\{
%E[
%+\sqrt{t}z)
%]_{|z=Z_V}\}\\
&&\qquad=\frac{\sqrt{1-t}}{d_i}\sum_{j=1}^m E\biggl[
\frac{\partial^2\varphi}{\partial x_i\,\partial x_j}\bigl(\sqrt
{1-t}(I_{d_1} (h_1),\ldots,I_{d_m} (h_m))
+\sqrt{t}Z_V\bigr)\\
&&\hspace*{166pt}\qquad{}\times \langle D[I_{d_i} (h_i)],D[I_{d_j} (h_j)]\rangle_\HH
\biggr].
\end{eqnarray*}
Hence, $\Psi'(t)$ equals
\begin{eqnarray*}
&&\frac1{2}
\sum_{i,j=1}^m
E \biggl[
\frac{\partial^2\varphi}{\partial x_i\,\partial x_j}\bigl( \sqrt
{1-t}(I_{d_1} (h_1),\ldots,I_{d_m} (h_m))
+ \sqrt{t}Z_V\bigr) \\
&&\hspace*{43pt}\qquad{}\times \biggl( V(i,j) -
\frac{1}{d_i}\langle D[I_{d_i} (h_i)],D[I_{d_j} (h_j)]\rangle_\HH
\biggr)
\biggr],
\end{eqnarray*}
so that we get
\begin{eqnarray*}
\delta_2&\leq& \|\varphi''\|_\infty
\sum_{i,j=1}^m
E\biggl[\biggl|V(i,j)-
\frac{1}{d_i}\langle D[I_{d_i} (h_i)],D[I_{d_j} (h_j)]\rangle_\HH
\biggr|\biggr]\\
&\leq& \|\varphi''\|_\infty
\sum_{i,j=1}^m
\sqrt{E\biggl[\biggl(V(i,j)-
\frac{1}{d_i}\langle D[I_{d_i} (h_i)],D[I_{d_j} (h_j)]\rangle_\HH
\biggr)^2\biggr]}\\
&=& \|\varphi''\|_\infty
\sum_{i,j=1}^m
\frac{1}{d_i}\sqrt{
{\operatorname{Var}}(
\langle D[I_{d_i} (h_i)],D[I_{d_j} (h_j)]\rangle_\HH)}.
\end{eqnarray*}

\textit{Step 4: Bounding ${\operatorname{Var}}(\langle D[I_{d_i}
(h_i)],D[I_{d_j} (h_j)]\rangle_\HH)$}.
Assume, for instance, that $i\leq j$. We have
\begin{eqnarray*}
&&\langle D[I_{d_i} (h_i)],D[I_{d_j} (h_j)]\rangle_\HH\\
&&\qquad=d_id_j\int_{0}^1
I_{d_i-1}(h_i(\cdot,a))I_{d_j-1}(h_j(\cdot,a))\,da\\
&&\qquad=d_id_j\int_{0}^1\sum_{r=0}^{d_i-1}r!\pmatrix{{d_i-1}\cr{r}}\pmatrix{{d_j-1}\cr{r}}
I_{d_i+d_j-2-2r}\bigl(h_i(\cdot,a)\,\widetilde{\otimes}_r\, h_j(\cdot
,a)\bigr)\,da\\
&&\qquad\hspace*{245pt}\mbox{(by Proposition \ref{P:MultiplictForm})}\\
&&\qquad=d_id_j\sum_{r=0}^{d_i-1}r!\pmatrix{{d_i-1}\cr{r}}\pmatrix{{d_j-1}\cr{r}}
I_{d_i+d_j-2-2r}(h_i\,\widetilde{\otimes}_{r+1}\, h_j) \\
&&\qquad=d_id_j\sum_{r=1}^{d_i}(r-1)!\pmatrix{{d_i-1}\cr{r-1}}\pmatrix{{d_j-1}\cr{r-1}}
I_{d_i+d_j-2r}(h_i\,\widetilde{\otimes}_{r}\, h_j).
\end{eqnarray*}
Hence, if $d_i<d_j$, then ${\operatorname{Var}}(\langle D[I_{d_i}
(h_i)],D[I_{d_j} (h_j)]\rangle_\HH)$ equals
\[
d_i^2d_j^2\sum_{r=1}^{d_i}(r-1)!^2\pmatrix{{d_i-1}\cr{r-1}}^2\pmatrix{{d_j-1}\cr{r-1}}^2
(d_i+d_j-2r)!\|
h_i\,\widetilde{\otimes}_{r}\, h_j\|^2_{\HH^{\otimes(d_i+d_j-2r)}},
\]
while, if $d_i=d_j$, it equals
\[
d_i^4\sum_{r=1}^{d_i-1}(r-1)!^2\pmatrix{{d_i-1}\cr{r-1}}^4
(2d_i-2r)!\|
\,h_i\widetilde{\otimes}_{r}\, h_j\|^2_{\HH^{\otimes(2d_i-2r)}}.
\]
Now, let us stress the two following estimates.
If $r<d_i\leq d_j$, then
\begin{eqnarray*}
\|h_i\,\widetilde{\otimes}_r\, h_j\|^2_{\HH^{\otimes(d_i+d_j-2r)}}
&\leq& \|h_i\otimes_r h_j\|^2_{\HH^{\otimes(d_i+d_j-2r)}}\\
&=&\langle h_i\otimes_{d_i-r} h_i, h_j\otimes_{d_j-r}h_j\rangle_{\HH
^{\otimes2r}}\\
&\leq&
\|h_i\otimes_{d_i-r}h_i\|_{\HH^{\otimes2r}}\|h_j\otimes_{d_j-r}h_j\|
_{\HH^{\otimes2r}}\\
&\leq&\tfrac12(
\|h_i\otimes_{d_i-r}h_i\|_{\HH^{\otimes2r}}^2+\|h_j\otimes
_{d_j-r}h_j\|
_{\HH^{\otimes2r}}^2
).
\end{eqnarray*}
If $r=d_i<d_j$, then
$
\|h_i\,\widetilde{\otimes}_{d_i}\, h_j\|^2_{\HH^{\otimes(d_j-d_i)}}
\leq
\|h_i\otimes_{d_i} h_j\|^2_{\HH^{\otimes(d_j-d_i)}} \leq
\|h_i\|^2_{\HH^{\otimes d_i}}\break\|h_j\otimes_{d_j-d_i}h_j\|_{\HH
^{\otimes2d_i}}.
$
By putting\vspace*{1pt} all these estimates in the previous expression for
${\operatorname{Var}}(\langle D[I_{d_i} (h_i)],D[I_{d_j} (h_j)]\rangle
_\HH)$, we get, using
also Lemma \ref{Lem:ContNorms}, that
$
\frac1{d_i}\sqrt{{\operatorname{Var}}(\langle D[I_{d_i} (h_i)],D[I_{d_j}
(h_j)]\rangle_\HH)} \leq\Delta_{ij},
$
for $\Delta_{ij}$ defined by (\ref{deltaij}).
This completes the proof of the theorem.
\end{pf}

We now translate the bound in Theorem \ref{multimoo2} into a bound for
indicators of convex sets.

\begin{cor} \label{multimoo2cor}
Let the notation and assumptions from Theorem \textup{\ref{multimoo2}} prevail.
We consider the class ${\mathcal H}(\mathbb{R}^m)$ of indicator
functions of
measurable convex sets
in $\mathbb{R}^m$.
Let
$
B_1= \frac12 \sum_{i=1}^m
\Delta_{ii}+\sum_{1\leq i<j\leq m}\Delta_{ij}$
and
\[
B_2 = C \Biggl(\beta+ \sqrt{\frac8\pi}\Biggr)
\Biggl[
\sum_{j=1}^m \bigl(16\sqrt2\beta\bigr)^{{(d_j-1)}/3}d_j!
\Biggr]^3 \sqrt{\max_{1\leq j\leq m}\max_{1\leq
i\leq N_j}\operatorname{Inf}_i(f_j)}.
\]

\begin{enumerate}
\item
Assume that the covariance matrix $V$ is the $m \times m$ identity
matrix $I_m$.
Then
\begin{eqnarray*}
&&\sup_{h \in{\mathcal H}(\mathbb{R}^m)}| E [h
(Q^1(\mathbf{X}),\ldots
,Q^m(\mathbf{X}))] -
E [h(Z_V)] |\\
&&\qquad\leq 8 (B_1 + B_2)^{1/4} m^{3/8}.
\end{eqnarray*}
\item\label{nonnegcase}
Assume that $V$ is of rank $k \leq m$, and let $\Lambda
=\operatorname{diag}(\lambda_1,
\ldots, \lambda_k)$ be
the diagonal matrix with the nonzero eigenvalues of $V$ on the
diagonal. Let $B$ be a
$m \times k$ column orthonormal matrix (i.e., $B^TB = I_k$ and
$BB^T=I_m$), such that $V = B \Lambda B^T$, and let
$
b= \max_{i,j} ( \Lambda^{-1/2} B^T )_{i,j}.
$
Then\vspace*{1pt}
$
| E [h(Q^1(\mathbf{X}),\ldots,\break Q^m(\mathbf{X})
)] -
E [h(Z_V)] | \leq8 (b^2 B_1 + b^3 B_2)^{1/4} m^{3/8}
$
for all $h \in{\mathcal H}(\mathbb{R}^m)$.
\end{enumerate}
\end{cor}
\begin{rem}\label{ivan34}
1. Notice that
$\sup_{z \in\mathbb{R}^m}
| P[ (Q^1(\mathbf{X}),\ldots,Q^m(\mathbf{X})) \leq z
] -\break P
[Z_V \leq z] |
\leq\sup_{h \in{\mathcal H}(\mathbb{R}^m)}
| E [h
(Q^1(\mathbf{X}),\ldots,Q^m(\mathbf{X})
)] -
E [h(Z_V)] |$.
Thus, Corollary \ref{multimoo2cor} immediately gives a bound for
Kolmogorov distance.\vspace*{-6pt}
\begin{longlist}[2.]
\item[2.] By using the bound for $\delta_2$ derived in the proof of Theorem
\ref{multimoo2} above,
and following
the same line of reasoning as in the proof of Corollary \ref{multimoo2cor},
we have, by keeping the notation of Theorem \ref{multimoo2}, that if
$\Delta_{ij}\to0$ for all $i,j=1,\ldots,m$ and
$\max_{1\leq j\leq m}\max_{1\leq i\leq
N_j}\operatorname{Inf}_i(f_j) \to0$,
then $
( Q_{d_1}(N_1,f_1,\mathbf{G}),\ldots,\break Q_{d_m}(N_m,f_m,\mathbf{G})
)\to
\mathscr{N}_m(0,V)$
as $N_1,\ldots,N_j\to\infty$,
\textit{in the Kolmogorov distance}.
\end{longlist}
\end{rem}

\begin{pf*}{Proof of Corollary \protect\ref{multimoo2cor}}
First assume that $V$ is the identity matrix.
We partially follow \cite{Rinott1996}, and let $\Phi$
denote the standard normal distribution in $\mathbb{R}^m$, and $\phi$
the corresponding
density function.
For $h \in{\mathcal H}(\mathbb{R}^m)$, define the smoothing
$
h_t(x) = \int_{\mathbb{R}^m} h ( \sqrt{t} y + \sqrt{1-t} x
) \Phi
(dy)$,
$0 < t < 1$.
The key result, found, for example, in \cite{Goetze1991}, Lemma~2.11,
is that,
for any probability \mbox{measure}~$Q$ on $\mathbb{R}^m$, for any $W \sim Q$
and $Z
\sim\Phi$, and for any $0 < t < 1$, we  have that $
{ \sup_{h \in{\mathcal H}(\mathbb{R}^m)} | E [h(W)] -  E [h(Z)] | }
\leq\frac43 [\sup_{h \in{\mathcal H}(\mathbb{R}^m)} |
E [h_t(W)] -
E[h_t(Z)] | +
2 \sqrt{m} \sqrt{t} ] .
$
Similarly as in \cite{Loh2007}, page 24, put
$
u(x, t, z) = (2 \pi t)^{-{m}/{2}}\times\break \exp(- \sum_{i=1}^m
\frac{(z_i - \sqrt{1-t} x_i)^2}{2t} ) ,
$
so that $
h_t(x) = \int_{\mathbb{R}^m} h ( z ) u(x, t, z) \,dz.
$
Observe that $u(x,t,z)$ is the density function of the Gaussian vector
$Y\sim\mathscr{N}(0,tI_m)$ taken in $z-\sqrt{1-t}x$.
%$
%h_t(x) = \int_{\mathbb{R}^m} h ( z ) u(x, t, z)\, dz.
%$
%By dominated convergence, we may differentiate under the integral and
%obtain,
%%for the second derivatives,
% + \frac{1-t}{t^2} \int_{\R^m} h ( z ) (z_i - \sqrt{1-t} x_i)^2 u(x,
%t, z) dz
%and, for $i \ne j$,
%$
%(z_j - \sqrt{1-t} x_j) u(x, t, z) dz.
%$
%For the third partial derivatives,
%& + \frac{(1-t)^\frac32}{t^3} \int_{\R^m} h ( z ) (z_i - \sqrt{1-t}
%x_i)^3 u(x, t, z) dz
%and, for $i \ne j$,
% u(x, t, z) dz\\
% && + \frac{(1-t)^\frac32}{t^3} \int_{\R^m} h ( z ) (z_i - \sqrt{1-t}
%x_i)^2 (z_j - \sqrt{1-t} x_j) u(x, t, z) dz,
%and, for $i,j,k$ all distinct,
%%\[
%= \frac{(1-t)^\frac32}{t^3} \int_{\R^m} h ( z ) (z_i - \sqrt{1-t} x_i)
%(z_j - \sqrt{1-t} x_j) (z_k - \sqrt{1-t} x_k) %u(x, t, z) dz.
Because $0 \leq h(z) \leq1$ for all $z\in\mathbb{R}^m$,
we may bound\vspace*{-3pt}
$
| \frac{\partial^2h_t}{\partial x_i^2}(x) | \leq
\frac{1-t}{t}
+\frac{1-t}{t^2} E[Y_i^2]= \frac{2(1-t)}{t}.
$
Similarly, for $i \ne j$,
$
| \frac{\partial^2h_t}{\partial x_i\, \partial x_j} (x) |
\leq
\frac{1-t}{t^2}E[|Y_i|]E[|Y_j|]=\frac{2(1-t)}{\pi t}.
$
Thus, we have $\| h''_t\|_\infty\leq1/t\leq1/t^{3/2}$.
Bounding the third derivatives in a similar fashion yields, for all
$i,j,k$ not necessarily
distinct, that
$| \frac{\partial^3h_t}{\partial x_i\, \partial x_j\, \partial x_k}
(x)|$
is less or equal than
\begin{eqnarray*}
&&\frac{(1-t)^{3/2}}{t^3} \max\bigl\{ 3 E [|Y_i|] t + E [| Y_i|^3];\\
&&\qquad \hspace*{50pt}E [| Y_j|] t + E[Y_i^2] E [| Y_j|];E[|Y_i|] E[| Y_j|] E[| Y_k|]
\bigr\},
\end{eqnarray*}
and so $\|h'''_t\|_\infty\leq1/t^{3/2}$.
With \cite{Goetze1991}, Lemma~2.11, and Theorem \ref{multimoo2}, this
gives that
\begin{eqnarray*}
&&\sup_{h \in{\mathcal H}(\mathbb{R}^m)}
| E [h(Q^1(\mathbf{X}),\ldots,Q^m(\mathbf{X})
)] -
E [h(Z_V)] |\\
&&\qquad\leq \frac43 \Bigl[\sup_{h \in{\mathcal H}(\mathbb{R}^m)}
| E [h_t(Q^1(\mathbf{X}),\ldots,Q^m(\mathbf{X})
)]
- E[h_t(Z_V)] | +
2 \sqrt{m} \sqrt{t} \Bigr] \\
&&\qquad\leq \frac83 \sqrt{m} \sqrt{t} + \frac43 (B_1 + B_2)
t^{-3/2} .
\end{eqnarray*}
This function is minimized for $t = \sqrt{ { 3 (B_1 + B_2) }/{(2
\sqrt{m}) }}$, yielding the first assertion.
%
%as required, that
%%\beas
%| E [h
%(Q^1(\mathbf{X}),\ldots,Q^m(\mathbf{X}))
%] -
% [h(Z_V)] |
%&\leq& 8 (B_1 + B_2)^\frac14 m^\frac38.
For Point~\ref{nonnegcase}, write $W=(Q^1(\mathbf{X}),\ldots
,Q^m(\mathbf{
X}))$ for simplicity.
For $h\in{\mathcal H}(\mathbb{R}^m)$, we have
\begin{eqnarray*}
&&E [h(W)] - E [h(Z_V)]\\
&&\qquad =
E[h(B\Lambda^{1/2}\times\Lambda^{-1/2}B^TW)]
- E[h(B\Lambda^{1/2}\times\Lambda^{-1/2}B^T Z_V)] .
\end{eqnarray*}
Put $g(x) =h ( B \Lambda^{1/2} x) $. Then, $g \in{\mathcal
H}(\mathbb{R}^k)$
and, thanks to \cite{Goetze1991}, Lemma 2.11, we can write
\begin{eqnarray*}
&& \sup_{h \in{\mathcal H}(\mathbb{R}^m)} | E [h(W)] - E [h(Z_V)] |\\
&&\qquad\leq
\sup_{g \in{\mathcal H}(\mathbb{R}^k)} | E [g(\Lambda^{-1/2}
B^T W)] - E
[g(\Lambda^{-1/2}B^TZ_V)] |
\\
&&\qquad\leq \frac43 \Bigl[\sup_{g \in{\mathcal H}(\mathbb
{R}^k)} | E [g_t(\Lambda
^{-1/2} B^T W )]
- E[g_t(\Lambda^{-1/2}B^TZ_V)] | +
2 \sqrt{k} \sqrt{t} \Bigr].
\end{eqnarray*}
We may bound the partial derivatives of $f_t(x) = g_t(\Lambda^{-
1/2} B^T x) $ using the
chain rule and the definition of $b$, to give that
$
\| f_t'' \|_\infty\leq b^2 t^{-3/2}$ and
$\| f_t''' \|_\infty\leq b^3 t^{-3/2}.
$
Using Theorem \ref{multimoo2} and minimizing the bound in $t$ as before
gives the assertion; the only changes are that $B_1$ gets multiplied by
$b^2$ and $B_2$ gets multiplied by $b^3$.
%\rightqed
\end{pf*}

%s7.2 ###
\subsection{More universality}

Here, we prove a slightly stronger version of Theorem~\ref
{T:UNIVmultiv} stated in
Section \ref{univers}. Precisely,
we add the two conditions (2)
and (3), making the criterion contained in Theorem \ref
{T:UNIVmultiv} more
effective for potential applications.

\begin{theorem}\label{T:UNIVmultivbis}
We let the notation of Theorem \ref{T:UNIVmultiv} prevail.
Then, as $n\rightarrow\infty$, the following four conditions \textup{
(1)}--\textup{(4)} are equivalent:
\textup{(1)} The vector $\{Q^j(n,\mathbf{G}) \dvtx j=1,\ldots,m\}$ converges in
law to
$\mathscr{N}_m(0,V)$;
\textup{(2)} for all $i,j=1,\ldots,m$, we have
$E[Q^i(n,\mathbf{G}) Q^j(n,\mathbf{G})]\to V(i,j)$ and
$E[Q^i(n,\mathbf{G})^4]\to3V(i,i)^2$ as $n\to\infty$;
\textup{(3)}~for all $i,j=1,\ldots,m$, we have
$E[Q^i(n,\mathbf{G}) Q^j(n,\mathbf{G})]\to V(i,j)$ and,
for all $1\leq i\leq m$ and $1\leq r \leq
d_i-1$, we have
$\|f_n^{(i)} \star_r f_n^{(i)}\|_{2d_i-2r} \rightarrow0$;
\textup{(4)} for every sequence $\mathbf{X } = \{X_i \dvtx  i\geq1\}$ of
independent centered random variables,
with unit variance and such that $\sup_i E|X_i|^{3} <\infty$, the vector
$\{Q^j(n,\mathbf{X}) \dvtx  j=1,\ldots,m\}$ converges in law to $\mathscr
{N}_m(0,V)$ for the Kolmogorov distance.
\end{theorem}

For the proof of Theorem \ref{T:UNIVmultivbis}, we need the
following result, which consists in a collection of
some of the findings contained in the papers by Peccati and Tudor
\cite{PT}. Strictly speaking, the original statements contained in
\cite
{PT} only deal with positive
definite covariance matrices: however, the extension to a nonnegative
matrix can be easily achieved by using
the same arguments as in Step 3 of the proof of Theorem \ref{multimoo2}.

\begin{theorem}\label{T:GIO+CIPRIAN}
Fix integers $m\geq1$ and $d_m\geq\cdots\geq
d_1\geq1$.
Let $V=\{V(i,j)\dvtx i,j=1,\ldots,m\}$ be a $m\times m$ nonnegative symmetric
matrix. For any $n\geq1$ and $i=1,\ldots,m$, let $I_{d_i}
(h_{i}^{(n)})$ belong to
the $d_i${th} Gaussian chaos $C_{d_i}$.
Assume that
$F^{(n)}=(F^{(n)}_1,\ldots,F^{(n)}_m):=(I_{d_1}(h_{1}^{(n)}),
\ldots,I_{d_m}(h_{m}^{(n)}))$, $n\geq1,$ is such that
$
\lim_{n\to\infty}
E[F_i^{(n)}F_j^{(n)}]=V(i,j)$, $1\leq i,j\leq m.
$\vspace*{-2pt}
Then, as $n\to\infty$, the following four assertions \textup{(i)}--\textup{(iv)} are
equivalent:
\textup{(i)} For every $1\leq i\leq m$, $F_i^{(n)}$ converges
in distribution
to a centered Gaussian random variable with variance $V(i,i)$;
\textup{(ii)}~for every $1\leq i\leq m$, $E
[(F_i^{(n)})^4
]\rightarrow3V(i,i)^2$;
\textup{(iii)} for every $1\leq i\leq m$ and every
$1\leq r \leq d_i-1$,
$\|h_{i}^{(n)} \otimes_r h_{i}^{(n)}\|_{\EuFrak{H}^{\otimes
(2d_i-2r)}} \rightarrow0$;
\textup{(iv)} the vector $F^{(n)}$ converges in distribution to the
$d$-dimensional Gaussian
vector $\mathscr{N}_m(0,V)$.
\end{theorem}

\begin{pf*}{Proof of Theorem \protect\ref{T:UNIVmultivbis}}
The equivalences $\mathrm{(1)}\Leftrightarrow\mathrm{(2)} \Leftrightarrow
\mathrm{(3)}$ only consist in a reformulation of
the previous Theorem \ref{T:GIO+CIPRIAN}, by taking into account
the first identity in
Lemma \ref{Lem:ContNorms}
and the fact that (since we suppose that
the sequence $E[Q^j(n,\mathbf{G})^2]$ of variances is bounded, so that an
hypercontractivity
argument can be applied),
if Point \textup{(1)} is verified, then $\lim_{n\to\infty}
E[F_i^{(n)}F_j^{(n)}]=V(i,j)$ for all $1\leq i,j\leq m$. On
the other
hand, it is completely obvious that
\textup{(4)} implies \textup{(1)}, since $\mathbf{G}$ is a particular case of
such an $\mathbf{X}$. So, it remains to prove the
implication $\mathrm{(1)},\mathrm{(2)},\mathrm{(3)}\Rightarrow\mathrm{(4)}$.
%Let $\varphi:\R\to\R$ be a differentiable function with $\|\varphi\|_
Let $Z_V=(Z^1_V,\ldots,Z^m_V)\sim\mathscr{N}_m(0,V)$. We have
\begin{eqnarray*}
&&\sup_{z\in\mathbb{R}^m}| P[ Q^1(n,\mathbf{X})\leq
z_1,\ldots,
Q^m(n,\mathbf{X})\leq z_m
] \\
&&\hspace*{47pt}\qquad{}- P[Z^1_V\leq z_1,\ldots,Z^m_V\leq z_m
]|\leq
\delta^{(a)}_n+\delta^{(b)}_n
\end{eqnarray*}
with
\begin{eqnarray*}
\delta^{(a)}_n &=& \sup_{z\in\mathbb{R}^m} |
P[
Q^1(n,\mathbf{X}) \leq z_1,\ldots,
Q^m(n,\mathbf{X}) \leq z_m
] \\
&&\hspace*{23pt}{}- P[ Q^1(n,\mathbf{G}) \leq z_1,\ldots,
Q^m(n,\mathbf{G}) \leq z_m
]
|,\\
\delta^{(b)}_n &=& \sup_{z\in\mathbb{R}^m} |
P[
Q^1(n,\mathbf{G}) \leq z_1,\ldots,
Q^m(n,\mathbf{G}) \leq z_m
] \\
&&\hspace*{69pt}{}- P[Z^1_V \leq z_1,\ldots,Z^m_V
\leq z_m
] |.
\end{eqnarray*}
By assumption \textup{(3)}, we have
that $\Delta_{ij}\to0$ for all $i,j=1,\ldots,m$ [with $\Delta_{ij}$
defined by (\ref{deltaij})].
Hence, Remark \ref{ivan34} (Point 2) implies that
$\delta^{(b)}_n\to0$.
By assumption~\textup{(3)} (for $r=d_j-1$) and (\ref{EQ:cruxcontr})--(\ref
{multiINF}),\vspace*{-1.5pt}
we get that\break $\max_{1\leq i\leq N_n^{(j)}}\operatorname{Inf}_i(f_n^{(j)})\to0$
as $n\to\infty$ for all $j=1,\ldots,m$.
Hence, Corollary \ref{multimoo2cor} implies that $\delta^{(a)}_n\to0$,
which completes the proof.
\end{pf*}

%s8 ###
\section{Some proofs based on Malliavin calculus and Stein's
method}\label{S:PrOOFS}
%s8.1 ###
\subsection{The language of Malliavin calculus}
Let $\mathbf{G} = \{G_i\dvtx  i\geq1\}$ be an i.i.d. sequence of Gaussian
random variables with zero mean and unit
variance. In what follows,
we will systematically use the definitions and notation introduced in
Section \ref{S:WienerC}.
In particular, we shall encode
the structure of random variables belonging to some Wiener chaos by
means of increasing (tensor) powers of a
fixed real separable Hilbert space $\HH$. We recall that the first
Wiener chaos of $\mathbf{G}$ is the $L^2$-closed Hilbert space of random
variables of the type $I_1(h)$, where
$h \in\HH$. We shall denote by $L^2(\mathbf{G})$ the space of all
$\mathbb{R}$-valued random elements $F$
that are measurable with respect to $\sigma\{\mathbf{G}\}$ and verify
$E[F^2]<\infty$.
Also, $L^2(\Omega;\HH)$ denotes the space of all $\HH$-valued random
elements $u$,
that are measurable with respect to $\sigma\{\mathbf{G}\}$ and verify the
relation $E[\|u\|_\HH^2]<\infty$.
%We stress (see again \cite[Ch. 1]{Nbook} or \cite{Janson}) that any
%random variable $F$ belonging
%to $L^2(\mathbf{G})$ admits the following chaotic expansion:
%F=E[F]+\sum_{d=1}^\infty I_d(h_d),
%where the series converges in $L^2(\mathbf{G})$ and
%the symmetric kernels $h_d\in\HH^{\odot d}$, $d\geq1$, are
%uniquely determined by $F$. For convenience we also put $I_0(f_0) =
%E[F]$.
For the rest of this section, we shall use standard notation and
results from Malliavin calculus:
the reader is referred to \cite{Nbook} for a detailed presentation of
these notions. In particular, $D^m$ denotes the $m$th Malliavin
derivative operator, whose domain is denoted by $\mathbb{D}^{m,2}$ (we
also write $D^1 =D$). An important property of $D$ is that it satisfies
the following \textit{chain
rule}: if $g\dvtx \mathbb{R}^n\rightarrow\mathbb{R}$ is continuously
differentiable
and has bounded partial derivatives, and if
$(F_1,\ldots,F_n)$ is a vector of elements of ${\mathbb{D}}^{1,2}$,
then $g(F_1,\ldots,F_n)\in{\mathbb{D}}^{1,2}$ and
$
Dg(F_1,\ldots,F_n)=\sum_{i=1}^n
\frac{\partial g}{\partial x_i} (F_1,\ldots, F_n)DF_i.
$
One can also show that the chain rule
continues to hold when $(F_1,\ldots, F_n)$
is a vector of multiple integrals (of possibly different orders) and
$g$ is a polynomial in $n$ variables.
We denote by $\delta$ the
adjoint of the operator $D$, also called the \textit{divergence
operator}. If a random element $u \in L^{2}(\Omega; \HH)$
belongs to the domain of $\delta$, noted $\operatorname{Dom}\delta$, then the
random variable $ \delta(u)$ is
defined by the duality relationship
$E (F \delta(u))= E \langle D F, u \rangle_{\HH}$,
which holds for every $F \in{\mathbb{D}}^{1,2}$.
As shown in \cite{NouPecptrf}, if $F=I_d(h)$, with $h\in\HH^{\odot
d}$, then one can deduce by integrating by parts (and by an appropriate
use of \textit{Ornstein--Uhlenbeck operators}) that, for every $G\in
\mathbb
{D}^{1,2}$ and every
continuously differentiable $g\dvtx \mathbb{R}\to\mathbb{R}$
with a bounded derivative, the following important relations hold:
%
%
%e33 ###
\begin{eqnarray}\label{Mexico2}
E[g(G)F]&=&\frac1d
E[g'(G)\langle DG, DF \rangle_\HH]\quad  \mbox{and}
\nonumber
\\[-8pt]
\\[-8pt]
\nonumber
 E[GF] &=&
\frac1d
E[\langle DG, DF \rangle_\HH].
\end{eqnarray}

Let $h\in\HH^{\odot d}$ with $d\geq2$, and let $s\geq0$
be an integer.
The following identity is obtained by taking $F=I_d(h)$
and $G=F^{s+1}$ in the second formula of (\ref{Mexico2}),
and then by applying the chain rule:
%
%
%e34 ###
\begin{equation}\label{MallMoments}
E[I_d(h)^{s+2}] = \frac{s+1}{d}E[I_d(h)^s\|DI_d(h)\|
_\HH
^2].
\end{equation}

%s8.2 ###
\subsection{Relations following from Stein's method}
Originally introduced in\break \cite{Steinorig,Steinbook}, \textit{Stein's
method} can be described as a collection of probabilistic techniques,
allowing to compute explicit bounds on the distance between the laws of
random variables by means of differential operators. The reader is
referred to %\cite{Chen_Shao} and
\cite{Reinertsur}, and the references therein, for an introduction to
these techniques.
The following statement contains four bounds which can be obtained by
means of a combination of Malliavin
calculus and Stein's method. Points 1, 2 and 4 have been proved in
\cite
{NouPecptrf}, whereas the content of
Point 3 is new. Our proof of such a bound gives an explicit example of
the interaction between Stein's method
and Malliavin calculus. We also introduce the following notation: for
every $F=I_d(h)$,
we set $
T_0(F) = \sqrt{\operatorname{Var}(\frac1d \| DF\|^2_{\HH})}.
$

\begin{prop}\label{P:BBBoou} Consider $F=I_d(h)$ with $d\geq1$ and
$h\in\HH^{\odot d}$, and let $Z$ and $Z_\nu$ have
respectively a $\mathscr{N}(0,1)$ and a $\chi^2(\nu)$ distribution
($\nu
\geq1$). We have the following:

\begin{enumerate}[1.]
\item[1.] If $E(F^2)=1$, then $d_{\mathrm{TV}}(F,Z) \leq2 T_0(F)$,
$d_{\mathrm{W}}(F,Z) \leq T_0(F)$ and,
for every thrice differentiable function $\varphi\dvtx \mathbb{R}\to
\mathbb{R}$ such that $\|
\varphi'''\| < \infty$,
$
|E[\varphi(F)] - E[\varphi(Z)]| \leq C_* \times T_0(F),
$
where $C_*$ is given in \textup{(\ref{Cost2})}.

\item[2.] If $E(F^2)=2\nu$, then
\[d_{\mathrm{BW}}(F,Z_\nu) \leq\max\Biggl\{
\sqrt{\frac{2\pi}{\nu}},\frac1\nu+\frac2{\nu^2}
\Biggr\} \sqrt{E\biggl[\biggl(2\nu+2F-\frac1d \| DF\|^2_{\HH}\biggr)^2\biggr]}.
\]
\end{enumerate}
\end{prop}

\begin{pf}
Point 2 is proved in \cite{NouPecptrf}, Theorem 3.11.
Point 1 is proved in \cite{NouPecptrf}, Theorem 3.1,
except the bound for $|E[\varphi(F)] - E[\varphi(Z)]|$.
To prove it,
fix $\varphi$ as in the statement, and consider the \textit{Stein equation}
$
f'(x) -xf(x) = \varphi(x) -E[\varphi(Z)]$, $x\in\mathbb{R}.
$
It is easily seen that a solution is given by
$
f(x)=f_{\varphi}(x) = e^{x^2/2} \int_{-\infty}^x (\varphi(y) -
E[\varphi
(Z)]) e^{-y^2/2}\, dy.
$
Set $K_*=C_*\times[4\sqrt{2}(1+5^{{3d}/{2}})]^{-1}$ with
$C_*$ given by (\ref{Cost2}).
According to the forthcoming Lemma \ref{Lemma:steinbound}, we have
$
|f_\varphi' (x)|\leq K_* (1+|x|+|x|^2+|x|^3).
$
Now use (\ref{Mexico2}) with $g=f_\varphi$ and $G=F$,
as well as a standard approximation argument to take into account that
$f'_\varphi$
is not necessarily bounded, in order to write
\begin{eqnarray*}
&& |E[\varphi(F)] - E[\varphi(Z)]| \\
&&\qquad= |E[f'_\varphi(F) - Ff_\varphi(F)]|\\
&&\qquad= \biggl|E\biggl[f'_\varphi(F)\biggl(1 - \frac1d\|DF\|^2_{\HH}
\biggr)
\biggr]\biggr|\\
&&\qquad\leq K_*E\biggl[
(1+|F|+|F|^2+|F|^3)\biggl|1 - \frac1d\|DF\|^2_{\HH}\biggr|
\biggr]\\
&&\qquad\leq4K_* E\biggl[(1+|F|^3)\biggl|1 - \frac1d\|DF\|
^2_{\HH}
\biggr|\biggr].
\end{eqnarray*}
By applying Cauchy--Schwarz, by using
$E[(1+|F|^3)^2]\leq2(1+E[F^6])$,
and finally by exploiting Proposition \ref{P:Hyper}, one infers
the desired conclusion:
\begin{eqnarray*}
4K_*E\biggl[(1+|F|^3)\biggl|1 - \frac1d\|DF\|^2_{\HH}
\biggr|
\biggr]
&\leq& C_* \sqrt{E\biggl[\biggl(1 - \frac1d\|DF\|^2_{\HH}
\biggr)^2\biggr]}\\
&=&C_*
T_0(F).
\end{eqnarray*}
\upqed\end{pf}

\begin{lemma}\label{Lemma:steinbound} The function $f_\varphi$
verifies
$|f_\varphi' (x)|\leq K_* (1+|x|+|x|^2+|x|^3)$.
%(\ref{creuzadema}).
\end{lemma}
\begin{pf}
We want to bound the quantity $\vert f_{\varphi}^{\prime}(
x) \vert$, where $\varphi$ is such that $\varphi
( x) =\varphi( 0) +\varphi^{\prime}(
0)
x+\varphi^{\prime\prime}( 0) x^{2}/2+R( x)$,
with $%
|R( x) |\leq\Vert\varphi^{\prime\prime
\prime
}\Vert_{\infty}\vert x\vert^{3}/6$.
Let $Z\sim\mathscr{N}(0,1)$.
We have
$f_{\varphi}^{\prime}( x) =A( x) +B(
x)$,
with $A( x) := \varphi( x) -E[\varphi(
Z)]$ and $B( x) := xf_{\varphi}( x)
$. It
will become clear\vspace*{1pt}
later on that our bounds on $\vert
f_{\varphi}^{\prime}( x) \vert$ do not depend on the
sign of $x$, so that in what follows we will only focus on the case
$x>0$%
. Due to the assumptions on $\varphi$, we have that
$
A( x) =\varphi^{\prime}( 0) x+\frac
{1}{2}\varphi
^{\prime\prime}( 0) x^{2}+R( x) +C
:= ax+bx^{2}+R( x) +C,
$
where
$-C =\frac{\varphi^{\prime\prime}( 0) }{2}+E[R(
Z)]$
[note that the term $\varphi( 0) $ simplifies]. Also, by
using\vspace*{-2pt} $E\vert Z\vert^{3}=\frac{2\sqrt{2}}{\sqrt{\pi
}}$ and
$E\vert Z\vert=\frac{\sqrt{2}}{\sqrt{\pi}}$, we obtain
$
\vert C\vert\leq\frac{\vert\varphi
^{\prime\prime
}( 0) \vert}{2}+\frac{\Vert\varphi
^{\prime
\prime
\prime}\Vert_{\infty}}{3}\frac{\sqrt{2}}{\sqrt{\pi}}:=
C^{\prime}
$
and (recall that $x>0$)
$
\vert A( x) \vert\leq\vert
\varphi^{\prime
}( 0) \vert x+\frac{1}{2}\vert\varphi
^{\prime
\prime
}( 0) \vert x^{2}+\frac{1}{6}\Vert\varphi
^{\prime
\prime\prime}\Vert_{\infty}x^{3}+C^{\prime}
=\vert a\vert x+\vert b\vert x^{2}+\gamma
x^{3}+C^{\prime}
$
with
$\gamma:= \frac{1}{6}\Vert\varphi^{\prime\prime\prime
}\Vert_{\infty}.
$
On the other hand, since $E[A( Z)] =0$ by construction,
$
\vert B( x) \vert=xe^{{x^{2}}/{2}}
\vert
\int_{x}^{+\infty}A( y) e^{-{y^{2}}/{2}}\,dy
\vert
\leq xe^{{x^{2}}/{2}}\int_{x}^{+\infty}[ C^{\prime
}+
\vert
a\vert y+\vert b\vert y^{2}+\gamma y^{3}]
e^{-{y^{2}}/{2}}\,dy
:= Y_{1}( x) +Y_{2}( x) +Y_{3}(
x) +Y_{4}( x).$ We now evaluate the four terms $Y_{i}$ separately (observe that each of them
is positive):%
\begin{eqnarray*}
Y_{1}( x) &=&C^{\prime}xe^{{x^{2}}/{2}}\int
_{x}^{+\infty
}e^{-{y^{2}}/{2}}\,dy\leq C^{\prime}e^{{x^{2}}/{2}}\int
_{x}^{+\infty
}ye^{-{y^{2}}/{2}}\,dy=C^{\prime}; \\
Y_{2}( x) &=&xe^{{x^{2}}/{2}}\int_{x}^{+\infty
}
\vert
a\vert ye^{-{y^{2}}/{2}}\,dy=\vert a\vert x; \\
Y_{3}( x) &=&xe^{{x^{2}}/{2}}\int_{x}^{+\infty
}
\vert
b\vert y^{2}e^{-{y^{2}}/{2}}\,dy=|b|\biggl( x^2+xe^{
{x^{2}}/{2}%
}\int_{x}^{+\infty} e^{-{y^{2}}/{2}}\,dy\biggr)\\
&\leq&\vert b\vert( x^{2}+1); \\
Y_{4}( x) &=&xe^{{x^{2}}/{2}}\int_{x}^{+\infty
}\gamma
y^{3}e^{-{y^{2}}/{2}}\,dy=\gamma x( x^{2}+2) =\gamma
x^{3}+2\gamma x.
\end{eqnarray*}

By combining the above bounds with $|f'_\varphi(x)|\leq|A(x)|+|B(x)|$,
one infers that%
\begin{eqnarray*}
\vert f_{\varphi}^{\prime}( x) \vert
&\leq
&2C^{\prime}+\vert b\vert+x( 2\vert a
\vert
+2\gamma) +x^{2}\vert b\vert+x^{3}2\gamma\\
&\leq&\max\{ 2C^{\prime}+\vert b\vert
;2\vert
a\vert+2\gamma;\vert b\vert;2\gamma\}
\times
( 1+x+x^{2}+x^{3}) \\
&=&\max\{ 2C^{\prime}+\vert b\vert;2\vert
a\vert+2\gamma\} \times( 1+x+x^{2}+x^{3}) ,
\end{eqnarray*}
which yields the desired conclusion.
\end{pf}

%s8.3 ###
\subsection{\texorpdfstring{Proof of Theorem \protect\ref{T:4thcumulant}}{Proof of Theorem 3.1}}\label{SS:Proof4thCUM}
Let $F=I_d(h)$, $h\in\HH^{\odot d}$. In view of Proposition~\ref
{P:BBBoou}, it is sufficient to show that
$
T_0(F) = T_1(F) \leq T_2(F).
$
Relation (3.42) in~\cite{NouPecptrf} yields that
%
%
%e35 ###
\begin{equation}\label{Scott1}
\frac1d\Vert DF \Vert_{\HH} ^{2}=E(F^{2})+d\sum_{r=1}^{d-1}(
r-1) !\pmatrix{{d-1}\cr{r-1}}^{2}I_{2d-2r}(
h\,\widetilde{\otimes}_{r}\,h),
\end{equation}
which, by taking the orthogonality of multiple integrals of different
orders into account, yields\vspace*{-1pt}
$\operatorname{Var}(\frac1d \| DF\|^2_{\HH})=d^{2}\sum
_{r=1}^{d-1}( r-1)!^{2} \left({{d-1}\atop{r-1}}\right)
^{4}(2d-2r) !\| h\,\widetilde{\otimes}_{r}\break h\| _{\mathfrak
{H}^{\otimes2( d-r) }}^{2},
$
and so $T_0(F) = T_1(F)$. From Proposition \ref{P:MultiplictForm},
we get
$F^2=\break \sum_{r=0}^{d} r!\left({{d}\atop{r}}\right)^2 I_{2d-2r}(h\,\widetilde{\otimes
}_r\, h).
$
To conclude the proof, we use (\ref{MallMoments}) with $s=2$,
combined with the previous identities, as well as the assumption that
$E ( F^2) =1$,
to get that
\begin{eqnarray*}
E[F^4]-3&=&\frac3d E(F^2\|DF\|^2_{\HH}) - 3(d!\|h\|^2_{\HH^{\otimes
d}})^2\\
&=&3d\sum_{r=1}^{d-1}r!( r-1) !\pmatrix{{d}\cr{r}}^{2}
\pmatrix{{d-1}\cr{r-1}}^{2}( 2d-2r) !\| h\,\widetilde{\otimes}%
_{r}\,h\| _{\HH^{\otimes2( d-r) }}^{2}.
\end{eqnarray*}
Hence,
$
\operatorname{Var}(\frac1d \| DF\|^2_{\HH}) \leq\frac
{d-1}{3d}[E(F^4)-3],
$
thus yielding $T_1(F) \leq T_2(F)$.

%s8.4 ###
\subsection{\texorpdfstring{Proof of Theorem \protect\ref{T:3and4}}{Proof of Theorem 3.6}}\label{SS:Proof4thCUMbis}
Let $F=I_d(h)$, $h\in\HH^{\odot d}$. In view of Proposition~\ref{P:BBBoou}
and since $L^{-1}F=-\frac1dF$, it is sufficient to show that
\[
\sqrt{E\biggl[\biggl( 2\nu+2F-\frac1d\|DF\|^2_\HH\biggr)^2
\biggr]} =
T_3(F) \leq T_4(F).
\]
By taking into account the orthogonality of multiple integrals of
different orders,
relation (\ref{Scott1}) yields
\begin{eqnarray*}
&&E\biggl[\biggl( 2\nu+2F-\frac1d\|DF\|^2_\HH\biggr)^2\biggr]\\
&&\qquad=
4d!\biggl\| h-\frac{d!^2}{4(d/2)!^3}h\,\widetilde{\otimes
}_{d/2}\,h\biggr\|
^2_{\HH^{\otimes d}} \\
&&\qquad\quad{}+ d^{2}
\mathop{\sum_{r=1,\ldots,d-1}}_{r\neq d/2}
( r-1) !^{2}\pmatrix{{d-1}\cr{r-1}}
^{4}( 2d-2r) !\Vert h\,\widetilde{\otimes
}_{r}\,h
\Vert_{
\mathfrak{H}^{\otimes( 2d-2r) }}^{2},
\end{eqnarray*}
and, consequently, $T_3(F) = \sqrt{E[( 2\nu+2F-\frac1d\|
DF\|
^2_\HH)^2]}$.
On the other hand, by combining (\ref{MallMoments})
(for $s=1$ and $s=2$) with $F^2=\sum_{r=0}^{d} r!{\left({d}\atop{r}\right)}^2
I_{2d-2r}(h\,\widetilde{\otimes}_r \,h)
$ [see the proof of Theorem \ref{T:4thcumulant}], we get,
still by
taking into account the orthogonality of multiple integrals of
different orders,
\begin{eqnarray*}
&&E[F^4]-12E[F^3]\\
&&\qquad=
12\nu^2-48\nu+24d!\biggl\| h-\frac{d!^2}{4(d/2)!^3}h\,\widetilde
{\otimes
}_{d/2}\,h\biggr\|^2_{\HH^{\otimes d}} \\
&&\qquad\quad{}+3d
\mathop{\sum_{r=1,\ldots,d-1}}_{r\neq d/2}
r!( r-1) !
\pmatrix{{d}\cr{r}}^2
\pmatrix{{d-1}\cr{r-1}}^{2}( 2d-2r) !\Vert h\,\widetilde
{\otimes
}_{r}\,h\Vert_{
\mathfrak{H}^{\otimes( 2d-2r) }}^{2}.
\end{eqnarray*}
It is now immediate to deduce that $T_3(F) \leq T_4(F)$.

\section*{Acknowledgments}
Part of this paper was written while the
three authors
were visiting the Institute for Mathematical Sciences of the National
University of Singapore, in the occasion of the program ``Progress in
Stein's Method''
(January 5--February 6, 2009). We heartily thank Andrew Barbour, Louis
Chen and Kwok Pui Choi for their kind hospitality and generous support.
We would also like to thank an anonymous Associate Editor and an
anonymous referee for helpful comments.

% imsref loaded by akundreckaite, 2010-05-21 08:51:28

%

\printaddresses

\end{document}